\documentclass[final,10pt]{smj}

\newcommand{\cA}{{\cal A}}
\newcommand{\cG}{{\cal G}}

\newcommand{\cC}{{\cal C}}
\newcommand{\cE}{{\cal E}}
\newcommand{\cL}{{\cal L}}
\newcommand{\cP}{{\cal P}}

\newcommand{\bX}{\mathbf{X}}
\newcommand{\bY}{\mathbf{Y}}
\newcommand{\bx}{\mathbf{x}}
\newcommand{\bz}{\mathbf{z}}
\newcommand{\be}{\mathbf{e}}
\newcommand{\bpi}{\mathbf {\pi}} 
\newcommand{\bPi}{\Pi}
\newcommand{\cip}{\mbox{\,$\perp\!\!\!\perp$\,}} 
\newcommand{\cd}{\,|\,} 
\newcommand{\mX}{\mathbb{X}}
\newcommand{\adj}{\mathop{\rm adj}}
\newcommand{\pa}{\mathop{\rm pa}}
\newcommand{\notcip}{\mathbin{\not\!\!\mathord{\perp}\!\!\!\mathord{\perp}}}

\usepackage{epstopdf}

\Author{David Edwards\Affil{1} and Smitha Ankinakatte\Affil{1}}
\AuthorRunning{}

\Affiliations{  \item    Centre for Quantitative Genetics and Genomics,          Aarhus University,        Denmark.}

\CorrAddress{David Edwards,
             Centre for Quantitative Genetics and Genomics,
             Department of Molecular Biology and Genetics,
             Science and Technology,
             Aarhus University,
             Blichers All\^{e} 20,
             8830 Tjele,
             Denmark}
\CorrEmail{David.Edwards@agrsci.dk}
\CorrPhone{(+45)\;87157965}
\CorrFax{(+45)\;87154994}

\Title{Context-specific graphical models for discrete longitudinal data}
\TitleRunning{Context-specific graphical models}

\Abstract{
\cite{Ron1998} introduced a rich family of models for discrete longitudinal data, called acyclic probabilistic finite automata. These may be represented as directed graphs that embody context-specific conditional independence relations. Here the approach is developed from a statistical perspective. It is shown how likelihood ratio tests may be constructed using standard contingency table methods, a model selection procedure that minimizes a penalized likelihood criterion is described, and a way to extend the models to incorporate covariates is proposed. The methods are applied to a small-scale data set. Finally, it is shown that the models generalize certain subclasses of conventional undirected and directed graphical models.
}


\begin{document}

\maketitle

\section{Introduction}
\label{sec:introduction}
\cite{Ron1998} described an approach to the analysis of discrete longitudinal data using \emph{acyclic probabilistic finite automata} (APFA)\footnote{We use the same acronym for both the singular and plural forms (automaton and automata).}. Automata are mathematical objects used in computer science, for example to represent formal languages and regular expressions, and in machine learning, for tasks such as speech recognition, natural language processing, and machine translation. APFA are a specific type of automata that are suited as models for complex discrete longitudinal data. They are widely used to process high-dimensional genomic data, as they underlie the popular program Beagle \citep{browning2007a}, but do not otherwise appear to have been taken up by the statistical community. \cite{smith2008} independently introduced a closely related class called chain event models, as we sketch in Section~\ref{sec:modelclasses}.

\begin{figure}[!ht]
\begin{center}
\includegraphics[trim=2cm 2cm 0cm 2cm, clip=true, width=3in]{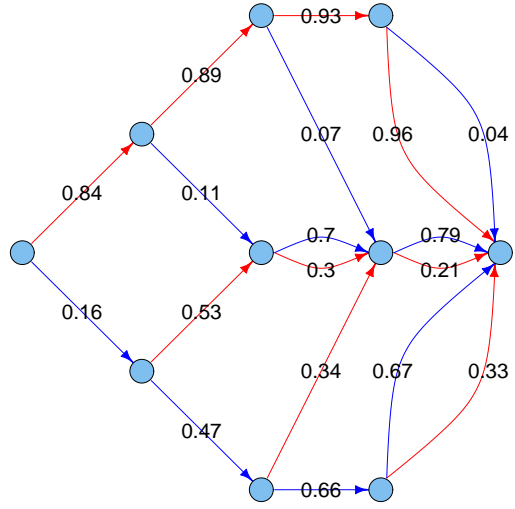}
\caption{An APFA for the wheezing data. The blue and red edges represents the presence and absence of wheezing, respectively, and the edge labels show transition probabilities. }
\label{fig:wheeze}
\end{center}
\end{figure}

To introduce the models we describe one for a small-scale longitudinal data set with four binary variables, denoted $W_1, \dots, W_{4}$, that record the presence or absence of wheezing at ages 7, 8, 9, 10 in a sample of 537 children (see Section~\ref{sec:application} below for more details). Figure~\ref{fig:wheeze} shows an APFA for this data set, selected using an algorithm we describe in Section~\ref{sec:modelselect} below. The variables are coded as 1 (absence) or 2 (presence): these values are represented in Figure~\ref{fig:wheeze} as red and blue edges, respectively. The edge labels show transition probabilities. Time flows from left to right, so we see, for example, that $\Pr(W_1= 2)=0.16$ and that $\Pr(W_2=2| W_1=2)=0.47$. We also see that
\begin{equation}
 \Pr(W_{3}=2|(W_1,W_2)) = \left\{ \begin{array}{ll}
                                        0.66 & \mbox{if} \ W_1=W_2=2 \\
                                        0.07 &  \mbox{if} \ W_1=W_2=1  \\
                                        0.70 & \mbox{otherwise}
                                      \end{array} \right.
\label{eq:condp1}
\end{equation}
and
\begin{equation}
 \Pr(W_{4}=2|(W_1,W_2,W_3)) = \left\{ \begin{array}{ll}
                                        0.67 & \mbox{if} \ W_1=W_2=W_3=2 \\
                                        0.04 &  \mbox{if} \ W_1=W_2=W_3=1  \\
                                        0.79 & \mbox{otherwise}
                                      \end{array} \right.
\label{eq:condp2}
\end{equation}
These properties can also be expressed in terms of conditional independence. For example (\ref{eq:condp2}) implies that $W_{4}$ is independent of $(W_2,W_2,W_3)$, given that $(W_1,W_2,W_3) \not \in \{(1,1,1),(2,2,2)\}$. This is an example of a context-specific conditional independence relation.

The structure of the paper is as follows. Section 2 introduces the model class. Sample trees are described in Section 3, and Section 4 treats maximum likelihood estimation. State merging is explained in Section 5, and related to the idea of model inclusion. Section 6 studies likelihood ratio tests between nested APFA: these turn out to be closely related to likelihood ratio tests in contingency tables. Section 7 exploits these results to modify the model selection algorithm of \cite{Ron1998} so as to optimize a penalized likelihood criterion. Section 8 indicates how the models may be extended to include covariates, and Section 9 applies the methods described to the wheezing data. In Section 10 APFA are related to conventional graphical models and to chain event graphs. The final section contains a discussion.

\section{Acyclic Probabilistic Finite Automata}

\begin{figure}[bht]
\begin{center}
\includegraphics[trim=2cm 0cm 0cm 2cm, clip=true,width=2.1in]{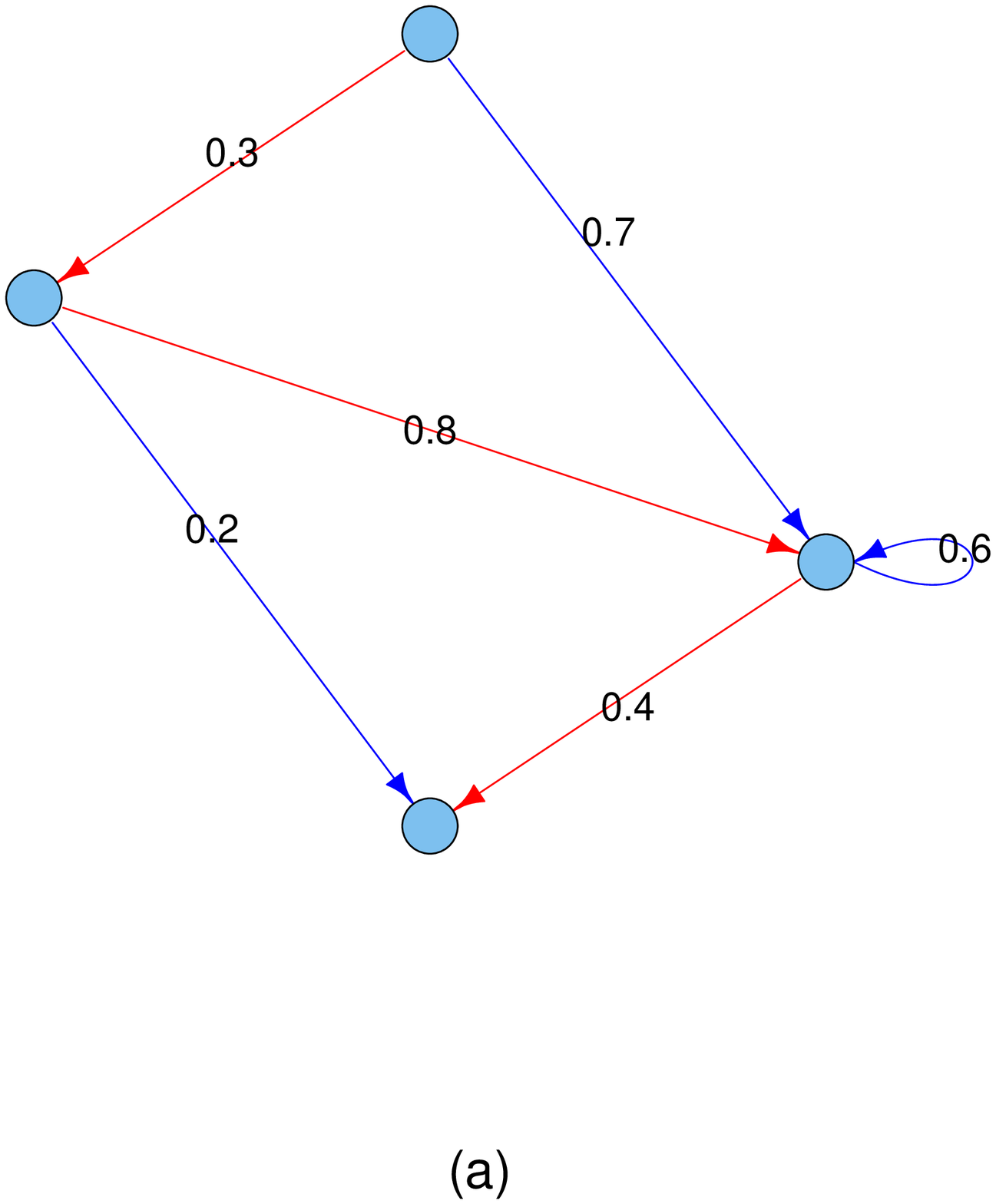}
\includegraphics[trim=2cm 0cm 0cm 2cm, clip=true,width=2.1in]{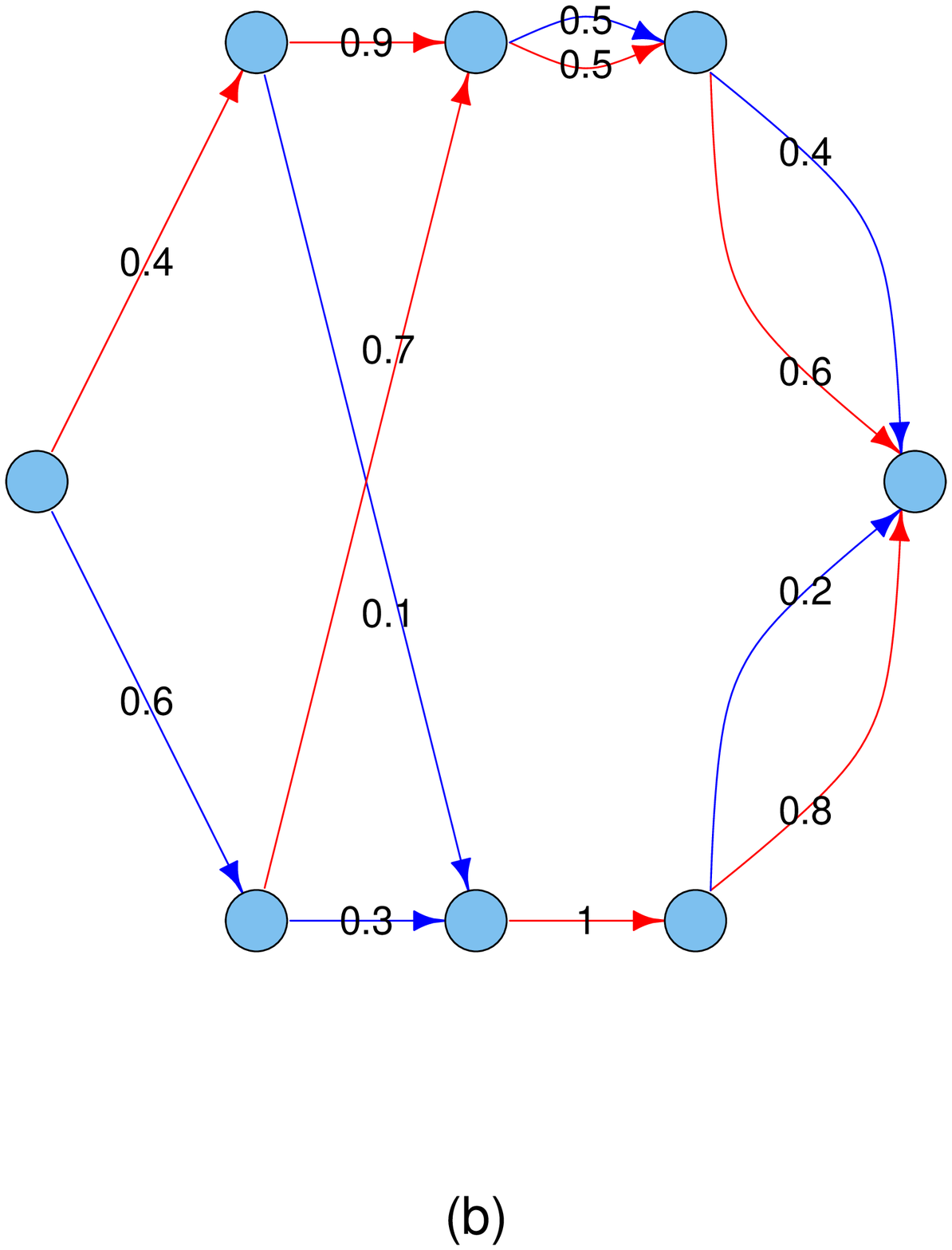}
\caption{Two multigraphs.  The first graph (a) is a probabilistic finite automaton (PFA). The second graph (b) is an acyclic probabilistic finite automaton (APFA). The colour of the edges denotes the symbol associated with the edge (red is "1" and blue is "2"). The edge labels are probabilities. }
\label{fig:simpleAPFA}
\end{center}
\end{figure}

We first describe the more general class of probabilistic finite automata (PFA). These are essentially devices that generate symbol strings. They are found in many variants, and the underlying theory has wide ramifications \citep{vidal2005a, vidal2005b}. A PFA may be represented as a directed multigraph, that is, a directed graph in which there may be multiple edges between node pairs. Figure~\ref{fig:simpleAPFA}a shows an example with four nodes or \emph{states}. The graph is required to have certain properties. It has precisely one initial or \emph{root} state with only outgoing edges, and precisely one final or \emph{sink} state with only incoming edges. All other states have at least one incoming edge and at least one outgoing edge. Self-loops (edges from a state to itself) are allowed. For every state there is a path from the root to the state, and a path from the state to the sink. Each edge $e$ has an associated symbol $\sigma(e)$ and a probability $\pi(e)$. The outgoing edges from each state have distinct symbols and the sum of their probabilities is one.

Such a graph defines a data generating process, which starts at the root, randomly chooses an outgoing edge according to their probabilities, generates the symbol associated with the edge, and traverses the edge to the next state: these steps are repeated until the process reaches the sink. In this way the graph defines a probability distribution over a set of output strings of possibly varying length. Each such string corresponds to a path from root to sink, and the probability of it being generated is the product of the edge probabilities along the path. Since there may be multiple edges between node pairs, edges cannot be uniquely identified as node pairs, nor paths as sequences of nodes. But since outgoing edges from each node have distinct symbols, edges are uniquely identified by their source node and symbol.

More formally, a PFA may be represented as a 7-tuple $(V,E,\Sigma,s,t,\sigma,\pi)$ where
$V$ is a set of vertices or nodes (often called states);
$E$ is a set of directed edges;
$\Sigma$ is a set of symbols (an alphabet);
$s$ and $t$ are maps $s: E \rightarrow V$ and $t: E \rightarrow V$ assigning to each edge its source and target nodes;
$\sigma$ is a map $s: E \rightarrow \Sigma$ assigning to each edge its symbol; and
$\pi$ is a map $\pi: E \rightarrow [0, 1]$ assigning to each edge its probability.

\label{sec:apfa}
PFA that generate strings of \emph{constant} length are termed APFA. These strings can be regarded as realizations of a discrete-valued random vector of fixed length. The defining graphical characterisation of an APFA is that all root-to-sink paths have the same length. This implies that all paths from the root to any specific state have the same length: this is called the \emph{level} of the state\footnote{Since acyclicity is necessary but not sufficient to ensure that all root-to-sink paths have the same length, the class could more aptly be called \emph{levelled} PFA.}. Each edge connects a state at one level to a state at the next level. Figure~\ref{fig:simpleAPFA}(b) shows an example with eight states and five levels (0 to 4).

Let $\cA$ be an APFA, let $p$ be the length of the root-to-sink paths in $\cA$, and let $\bX=(X_1, X_2, \dots X_p)$ be a vector of discrete random variables that take values in the sample spaces $\mX_i$, for $i=1,2 \dots p$.  Given a root-to-sink path $\be = (e_1, e_2, \dots ,e_p)$ in $\cA$ we equate the associated $p$-vector of symbols $\sigma(\be)= (\sigma(e_1), \sigma(e_2), \dots ,\sigma(e_p))$ with a realization of $\bX$. Distinct root-to-sink paths generate distinct symbol strings and hence distinct realizations of $\bX$.

The sample space of $\bX$ is given by $\mX(\cA)=\{\sigma(\be): \be \in \cE(\cA)\}$, where $\cE(\cA)$ is the set of root-to-sink paths in $\cA$.
Here $\mX(\cA)$ is some subspace of the product space $\mX = \prod \mX_i$.
For any $\bx \in \mX(\cA)$ we can find the unique root-to-sink path $\be$ such that $\bx=\sigma(\be)$: we write this as $\be=\sigma^{-1}(\bx)$. The sample space $\mX_i$ corresponds to the set of symbols generated by incoming edges to a level $i$ state.

The parameter vector $\bpi=\{\pi(e): e \in E(\cA)\}$ and the parameter space is
\[ \bPi = \{ \bpi: \pi(e) \geq 0 \ \forall e \in E(\cA) \ \mbox{and} \ \sum_{e : s(e)=v} \pi(e) =1 \ \forall v \in V(\cA)\}.  \]

The edge probabilities specify the marginal and conditional probabilities appearing in the standard factorization of the joint density of $\bX$
\begin{equation}
\Pr(\mathbf{X=x}) = \Pr(X_1=x_1) \prod_{i=2 \dots p} \Pr(X_i=x_i|X_{<i}=x_{<i})
\label{eq:factor}
\end{equation}
\noindent
where here and throughout we use shorthand expressions such as $\bX_{<i}=(X_1, \dots, X_{i-1})$, $\bx_{\geq i}=(x_{i}, \dots, x_p)$,
$\bY_{\geq i; \leq j}=(Y_i, \dots ,Y_j)$ and so forth.

When the data generating process arrives at state $w$ at level $i$, the distribution of the future observations $\bX_{>i}$ does not depend on the path the process took to arrive at $w$. This implies constraints on the joint distribution of $\bX$ which can be written as
\begin{equation}
\bX_{>i} \cip \bX_{\leq i} | \bX_{\leq i} \in {\cC(w)}
\label{eq:condind}
\end{equation}
where ${\cC(w)}=\{\sigma(\be): \be \in \cP(w)\}$, and $\cP(w)$ is the set of paths from the root to $w$.  For example, Figure~\ref{fig:simpleAPFA}(b) implies that
\[
(X_3,X_4) \cip X_{\leq 2}|X_{\leq 2} \in \{(1,2),(2,2)\}.
\]
For $i=1 \dots p$ let $I_i$ be a discrete random variable indicating which level $i$ node the root-to-sink path passes through. Then (\ref{eq:condind}) can be written more elegantly as
\begin{equation}
\bX_{>i} \cip \bX_{\leq i} | I_i= w.
\label{eq:condinda}
\end{equation}
\noindent
This is true for all level $i$ nodes $w$. So
\begin{equation}
\bX_{>i} \cip \bX_{\leq i} | I_i
\label{eq:condindb}
\end{equation}
\noindent
for $i=1, \dots, p-1$.  We can think of $I_i$ as representing the memory of the process at time $i$. For a more general form of (\ref{eq:condindb}) for chain event graphs, see \cite{thwaites2011}.

Thus an APFA expresses a set of conditional independence constraints on the distribution of $\bX$, and in this respect it resembles the dependence graph of a traditional graphical model \citep{Lauritzen1996, edwards2000}.  We compare the model classes in more detail in Section~\ref{sec:modelclasses}.

\begin{figure}[!ht]
\centering
\includegraphics[trim=2cm 2cm 1cm 2cm, clip=false, width=2.5in]{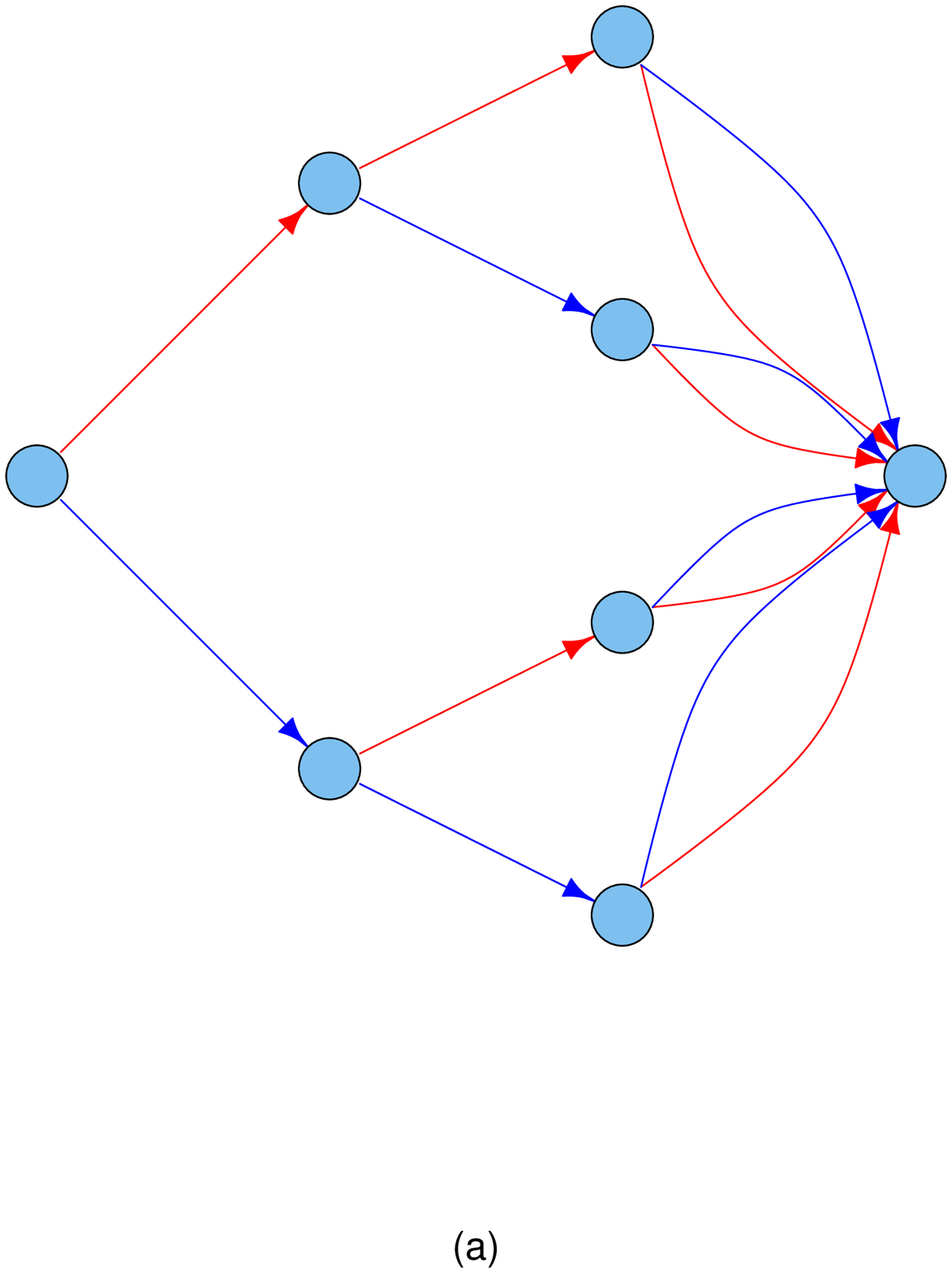}
\includegraphics[trim=0cm -7.4cm 0cm 7cm, clip=false, width=2.5in]{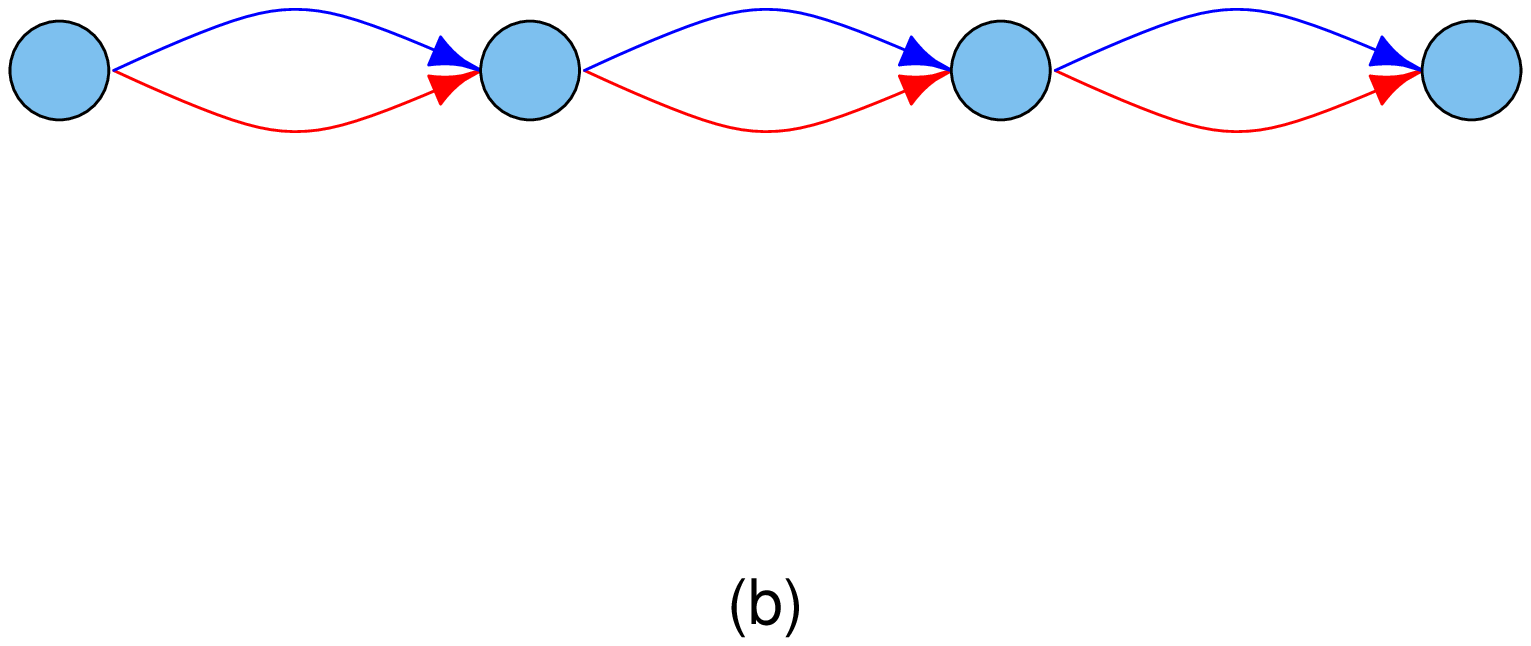}
\caption{The maximal, unrestricted APFA (a) and the minimal APFA (b) for three binary variables.}
\label{fig:maxmin}
\end{figure}

Absent edges in APFA may represent another type of constraint. For example, since state E has no outgoing blue edge,  Figure~\ref{fig:simpleAPFA}(b) implies that $\Pr(X_3=2|X_2=2)=0$. Here we should distinguish between \emph{structural} and \emph{random} zeroes. In some contexts it may be known prior to the analysis that when $X_2=2$, $X_3$ cannot be 2, so the absent edge represents a structural constraint built into the model. More commonly, perhaps, absent edges arise because the corresponding observations did not occur in the sample, and the APFA represents an estimate under a larger underlying model without the constraint. We make this idea more precise in Section~\ref{sec:merging} below.  In this paper we assume that absent edges reflect random rather than structural zeroes. We call an APFA $\cA$ \emph{complete}  if it has no absent edges, that is to say, each level $i$ node (except the sink) has $|\mX_{i+1}|$ outgoing edges. Otherwise we call it \emph{incomplete}. Clearly $\cA$ is complete if and only if $\mX(\cA) = \prod \mX_i$.

Figure~\ref{fig:maxmin} displays maximal and minimal APFA for three binary variables. The maximal model entails no restrictions on the joint distribution, and the minimal model represents complete independence.  Both are complete.

\section{Sample Trees}

\begin{figure}[!th]
\begin{center}
\includegraphics[trim=2cm 2cm 0cm 2cm, clip=true, width=4.5in]{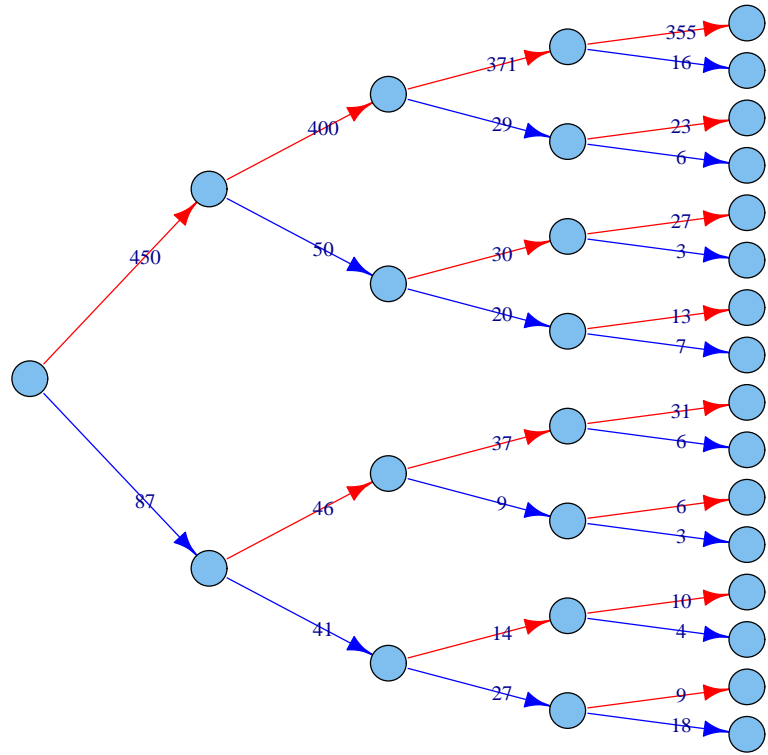}
\caption{The sample tree for the wheezing data. Red and blue edges correspond to the absence and presence of wheezing,
respectively.}
\label{fig:sampletree}
\end{center}
\end{figure}

Suppose that a data sample of the form $\bx^{(k)}=(x^{(k)}_1, \dots x^{(k)}_p)$ for $k=1 \dots N$ is available. The sample tree is a rooted tree in which the states represent partial outcomes $(\bx_{\leq q})$ for $q \in \{0, \dots, p\}$. The root of the tree represents the null outcome, and the nodes adjacent to the root the outcomes $x_1=1$, $x_1=2$ and so on. Each root-to-leaf path represents a distinct outcome $(x_1, \dots x_p)$ present in the sample data. Figure~\ref{fig:sampletree} shows the sample tree of the wheezing data set. The edges are labelled with the corresponding sample counts. Sample trees are called \emph{prefix tree acceptors} in machine learning and are closely related to \emph{tries} in computer science, and \emph{event trees} in Bayesian decision theory \citep{smith2008}. They provide a useful summary of discrete longitudinal data of small dimension.

If the leaves of the sample tree are contracted to a single node (the sink) an APFA is obtained: we call this the \emph{sample APFA}. It is typically used as start model in the selection algorithm described below in Section~\ref{sec:modelselect}.

Note that the sample APFA embodies no constraints of type (\ref{eq:condind}), but any states corresponding to partial outcomes $(\bx_{\leq q})$ not occurring in the data will be absent. So it is generally incomplete. It can be regarded as an estimate of the joint distribution under the unrestricted model, that is, the corresponding complete APFA in which each level $i$ node (except the sink) has $|\mX_{i+1}|$ outgoing edges.

Consider a data set with $N=1000$ observations of $p=100$ binary variables. The sample APFA can have at most $1000 \cong 2^{10}$ nodes at level $p-1$, but the corresponding complete APFA will have $2^{99}$. It clearly makes computational sense to exclude edges with zero counts from the sample APFA.

\section{Maximum Likelihood Estimation}
Suppose now that $\cA$ is an APFA such as that in Figure~\ref{fig:simpleAPFA}(b) whose edge probabilities $\pi(e)$ are unknown, and that independent samples $\bx^{(k)}=(x^{(k)}_1, \dots x^{(k)}_p)$ for $k=1 \dots N$ are drawn from $\cA$. We wish to estimate the $\pi(e)$. For $\bx \in \mX(\cA)$
\[
\Pr(\bx\cd \bpi) = \prod_{i = 1 \dots p} \pi(e_i)
\]
where $ \be= \sigma^{-1}(\bx)$, so that the likelihood of the sample is
\begin{equation}
 \prod_{k=1 \dots N} \Pr(\bx^{(k)}\cd \bpi) = \prod_{k=1 \dots N} \prod_{i = 1 \dots p} \pi(e_i^{(k)})
 \label{eq: likelihood}
\end{equation}
where  $ \be^{(k)}= \sigma^{-1}(\bx^{(k)})$. This can be re-written as
\[
 \prod_{k=1 \dots N} \Pr(\bx^{(k)}\cd \bpi) = \prod_{e \in E(\cA)}\pi(e)^{n(e)}
\]
where $n(e)$ is the \emph{edge count}, i.e. the number of observations in the sample whose root-to-sink path traverses the edge $e$.
We similarly define the \emph{node counts} $n(v)$ to be the number of observations in the sample whose root-to-sink path passes through  $v \in V$.

Maximum likelihood estimation is very straightforward: the edge probabilities are simply estimated as the relative frequencies of the corresponding counts. Thus for each $e \in E(\cA)$,
\begin{equation}
\hat{\pi}(e) =  \frac{n(e)}{n(v)},
\label{eq:edgeprob}
\end{equation}
where $v=s(e)$, the source node of $e$. The maximized log-likelihood under $\cA$ is
\begin{equation}
\hat{\ell}(\cA)=\sum_{e \in E(\cA)} n(e) \log \hat{\pi}(e).
\label{eq:loglike}
\end{equation}

Let $p(v)$ and $p(e)$ be the marginal probabilities of passing through a node $v \in V$ and an edge $e \in E(\cA)$, respectively. The maximum likelihood estimates of these quantities are the sample proportions, that is, $\hat{p}(v) = n(v)/N$ and $\hat{p}(e) = n(e)/N$.

\section{State Merging}
\label{sec:merging}

\begin{figure}[ht]
\centering
\includegraphics[trim=2cm 0cm 0cm 2cm, clip=true, width=2.5in]{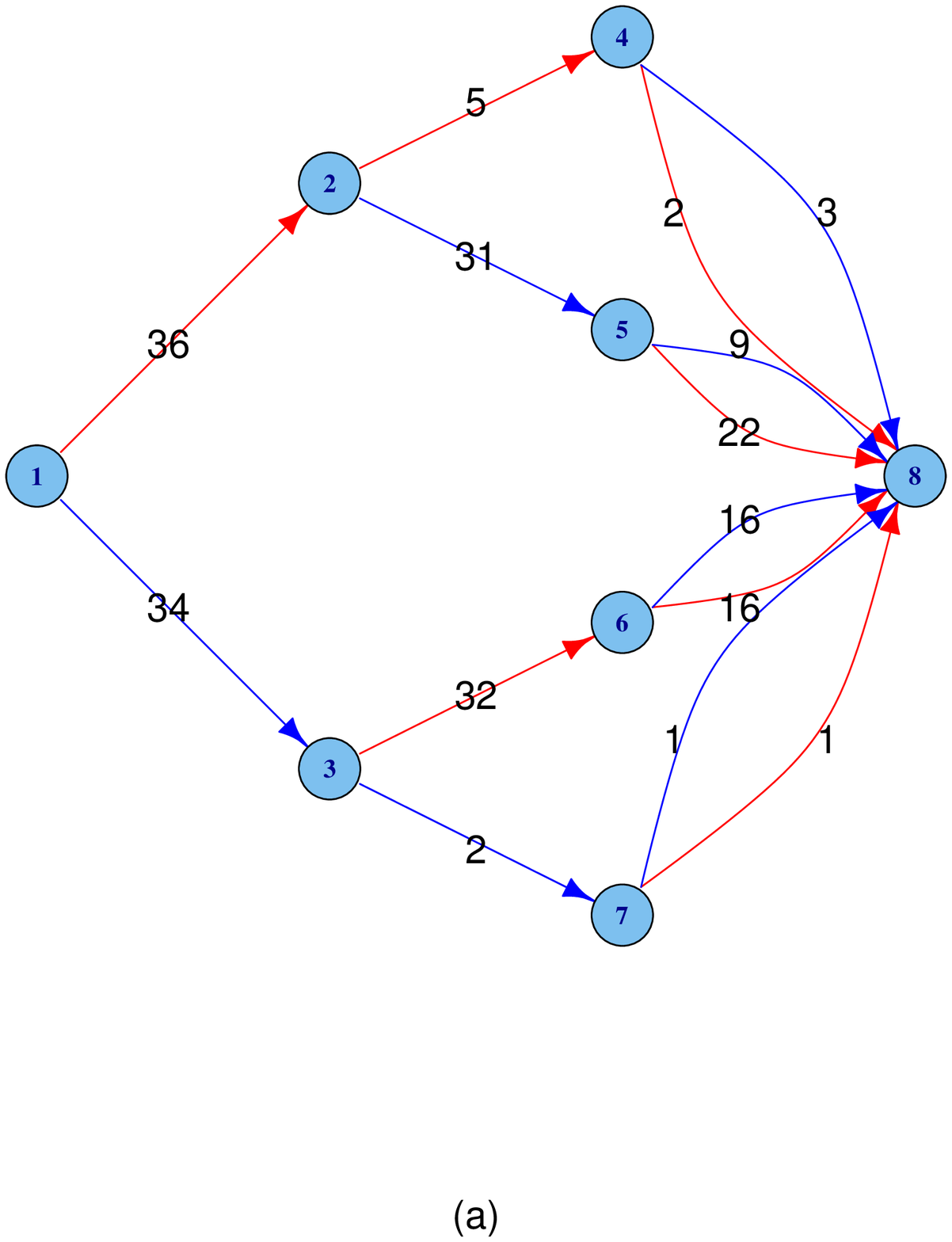}
\includegraphics[trim=2cm 0cm 0cm 2cm, clip=true, width=2.5in]{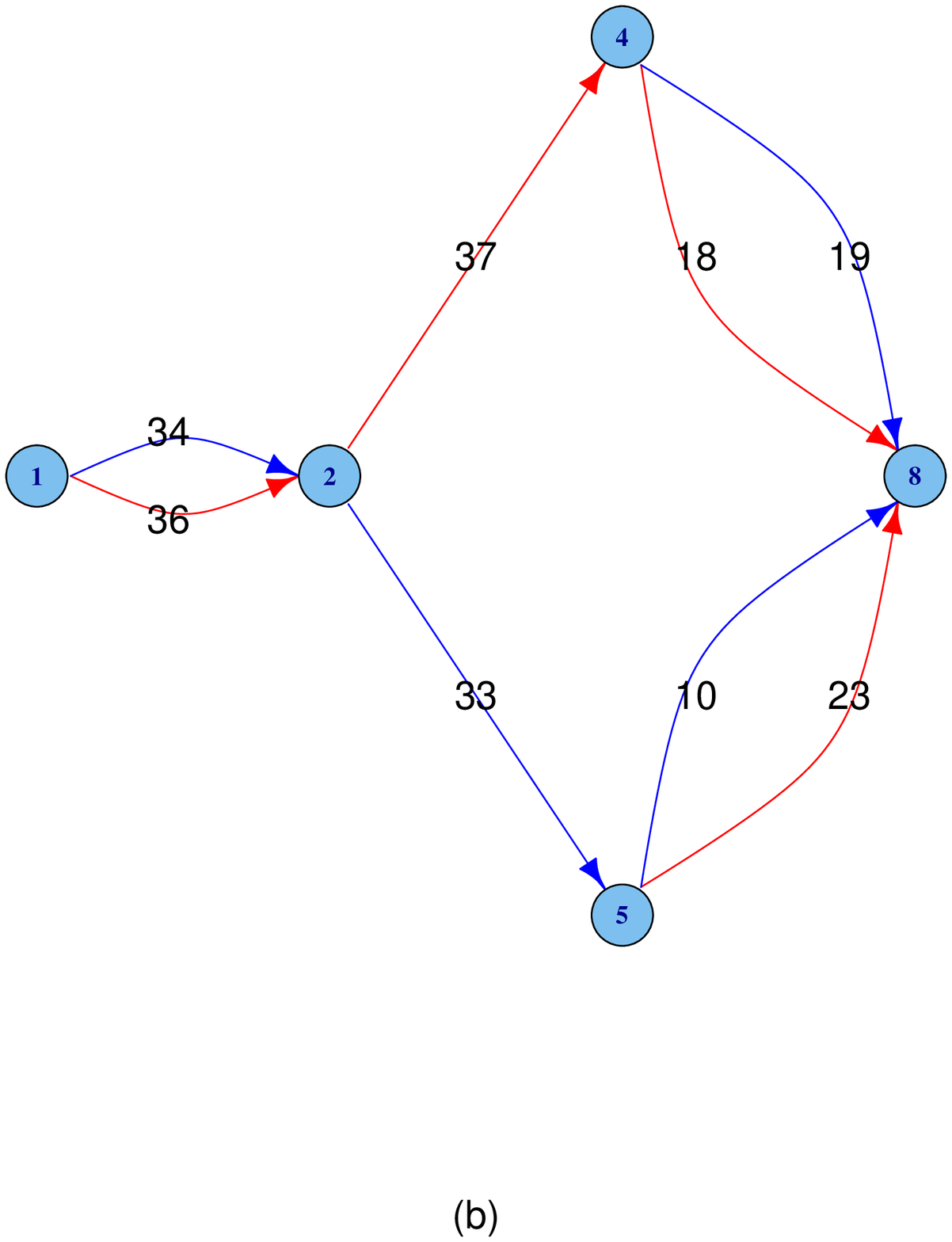}
\caption{An example of state merging. (a) shows an unrestricted APFA, and (b) shows the APFA obtained from (a) by merging node 3 with node 2,
which entails merging state 7 with state 5, and state 6 with state 4. The edge counts are also shown.}
\label{fig:merging}
\end{figure}

Simplifying parametric models typically involves setting parameters to be equal or zero. In contrast, simplifying APFA involves merging states. To retain the level structure, only states at the same level may be merged. Suppose we wish to merge state $w$ into state $v$: that is, redirect all incoming edges to $w$ to $v$, and redefine all outgoing edges from $w$ to outgo from $v$ instead. The former is unproblematic, but the latter may lead to the existence of outgoing edges from $v$ with duplicate symbols. Such edges must therefore also be merged, and if their target nodes are distinct, these must also be merged. The operation is thus recursive. For example, Figure~\ref{fig:merging}(b) is obtained from Figure~\ref{fig:merging}(a) by merging state 3 with state 2, which entails merging state 7 with state 5, and state 6 with state 4.

Using the same logic, more than two nodes may be merged. Let $s$ be a node set to be merged, and let $\cL(s)$ be the associated merge-list, that is, a list
containing $s$ and the other node sets that are merged. So for example the merge-list associated with merging $s=\{2,3\}$ in Figure~\ref{fig:merging}(a)  is $\cL(s) = \{2,3\}, \{5,7\}, \{4,6\}$.

A non-recursive characterization of state merging goes like this. Let $s$ be the node set to be merged. Then two nodes $x$ and $y$ are merged if and only if they are corresponding descendent nodes of two nodes $v$ and $w$ in $s$, that is, there exist paths $v \rightarrow x$ and $w \rightarrow y$ with identical symbol sequences. Similarly, two edges $e$ and $f$ are merged if and only if they are corresponding descendent edges of two nodes $v$ and $w$ in $s$, that is, there exist paths $v \rightarrow t(e)$ and $w \rightarrow t(f)$ with identical symbol sequences, whose last edges are $e$ and $f$.

Define the descendent subgraph of a node in an APFA to be the subgraph induced by the node and its descendants. Let $\cA$ and $\cA_0$ be the APFA shown in Figure~\ref{fig:merging}(a) and (b). Note that the descendent subgraphs of nodes 2 and 3 in $\cA$ and node 2 in $\cA_0$ are complete in the sense given in Section~\ref{sec:apfa}. This reflects that the conditional distributions of $(X_2,X_3)$ given $X_1=1$ and $2$ are unrestricted in $\cA$ and constrained to be equal in $\cA_0$. Thus $\cA_0$ is a submodel of $\cA$.

Here $\cA$ was complete. The incomplete case is more subtle, and is illustrated in Figure~\ref{fig:underspec}.
Let $\cA$ be the APFA shown in Figure~\ref{fig:underspec}(a). Nodes 2 and 3 in $\cA$ are merged to obtain $\cA_0$, shown in Figure~\ref{fig:underspec}(b). In $\cA$ node 3 has no outgoing blue edge, so under $\cA$, $\Pr(X_2=2|X_1=2)=0$, but under $\cA_0$ this does not hold. If the missing edge in $\cA$ represents a structural constraint, this constraint is not respected in $\cA_0$, so the merging is in conflict with this constraint. Moreover $\cA_0$ is not a submodel of $\cA$. It makes more sense to suppose that the missing edge represents a random zero, and that the model underlying $\cA$ is that shown in Figure~\ref{fig:merging}(a) (with different edge counts). Note that $\cA_0$ is a submodel of this.

To generalize this, define the \emph{completion} $\cA^+$ of an incomplete APFA $\cA$ as follows. Call a node $x$ of $\cA$ at level $i<p$ incomplete if it has fewer than $|X_{i+1}|$ outgoing edges. Recursively complete all incomplete nodes by adding the required number of new edges with the appropriate symbols: when $i<p-1$ this will require that the same number of new nodes are also added. In $\cA^+$ the descendent subgraphs of the nodes introduced are complete trees (with the final level collapsed to the sink). The following result is shown in Appendix A (in the online supplement).

\textbf{Assertion I: }\textit{When $\cA_0$ is obtained from an APFA  $\cA$ by state merging, the corresponding completed models are nested, that is to say, $(\cA_0)^+$ is a submodel of  $\cA^+$.}

Note that since completing an APFA only involves adding extra edges with zero edge counts, (\ref{eq:loglike}) implies that $\hat{\ell}(\cA) = \hat{\ell}(\cA^+)$ and $\hat{\ell}(\cA_0) = \hat{\ell}((\cA_0)^+)$.

\section{Hypothesis Testing}
\label{sec:testing}
These results can be used to construct likelihood ratio tests of nested hypotheses, that is of $\cA_0$ versus $\cA$, where $\cA_0$ is a submodel of $\cA$. For example, suppose that a sample of $N=70$ observations of $p=3$ binary variables is available, and consider the two APFA for these data shown in Figure~\ref{fig:merging}.

The likelihood ratio test (LRT) of $\cA_0$ versus $\cA$, often called the \emph{deviance}, is minus twice the logarithm of the likelihood ratio of $\cA_0$ versus $\cA$, that is
\begin{equation}
G^2 = - 2[\hat{\ell}(\cA) - \hat{\ell}(\cA_0)].
\label{eq:deviance}
\end{equation}
Here $\hat{\ell}(\cA)=-116.2117$ and $\hat{\ell}(\cA_0)=-142.7731$, so $G^2=53.1228$. Under $\cA_0$, $G^2$ is asymptotically $\chi^2(k)$ distributed where the degrees of freedom $k$ is the difference in model dimension (number of free parameters) between $\cA$ and $\cA_0$. This provides a goodness-of-fit test for the smaller model. By inspection of Figure~\ref{fig:merging} it appears that $\cA$ has 7 free parameters and $\cA_0$ has 4, so $k=3$, and clearly $\cA_0$ fits very poorly.

\begin{table}
\caption{\label{ctab1}Comparing the distribution of $(X_2,X_3)$ given nodes 2 and 3: data in Figure~\ref{fig:merging}(a).}
\centering
\fbox{%
\begin{tabular}{c|cccc}
  \hline
  source & (1,1) & (1,2) & (2,1) & (2,2) \\ \hline
  2 & 2 & 3 & 9 & 22 \\
  3 & 16 & 16 & 1 & 1 \\
  \hline
\end{tabular}}
\end{table}

The same test can be computed by applying a standard contingency table test of independence to the data in Table~\ref{ctab1}, which compares  the conditional distribution of $(X_2,X_3)$ given node 2 with that given node 3. The counts in the table are those of incoming edges to the sink in $\cA$. We may recall that for an $r \times c$ table of counts $\{n_{ij}\}_{i=1 \dots r; j=1 \dots c}$ the likelihood ratio test can be written as
\begin{equation}
G^2 = 2 \sum_{i,j} n_{ij} \log \frac{n_{ij} n_{++}}{n_{i+} n_{+j}}
\label{eq:G2}
\end{equation}
with degrees of freedom given as
\begin{equation}
 k = (\#\{i: n_{i+}>0\}-1)(\#\{j: n_{+j}>0\}-1)
\label{eq:df}
\end{equation}
where $n_{i+}$ and $n_{+j}$ are the row and column totals, respectively.

\begin{figure}[!ht]
\centering
\includegraphics[trim=2cm 0cm 0cm 2cm, clip=true, width=2.5in]{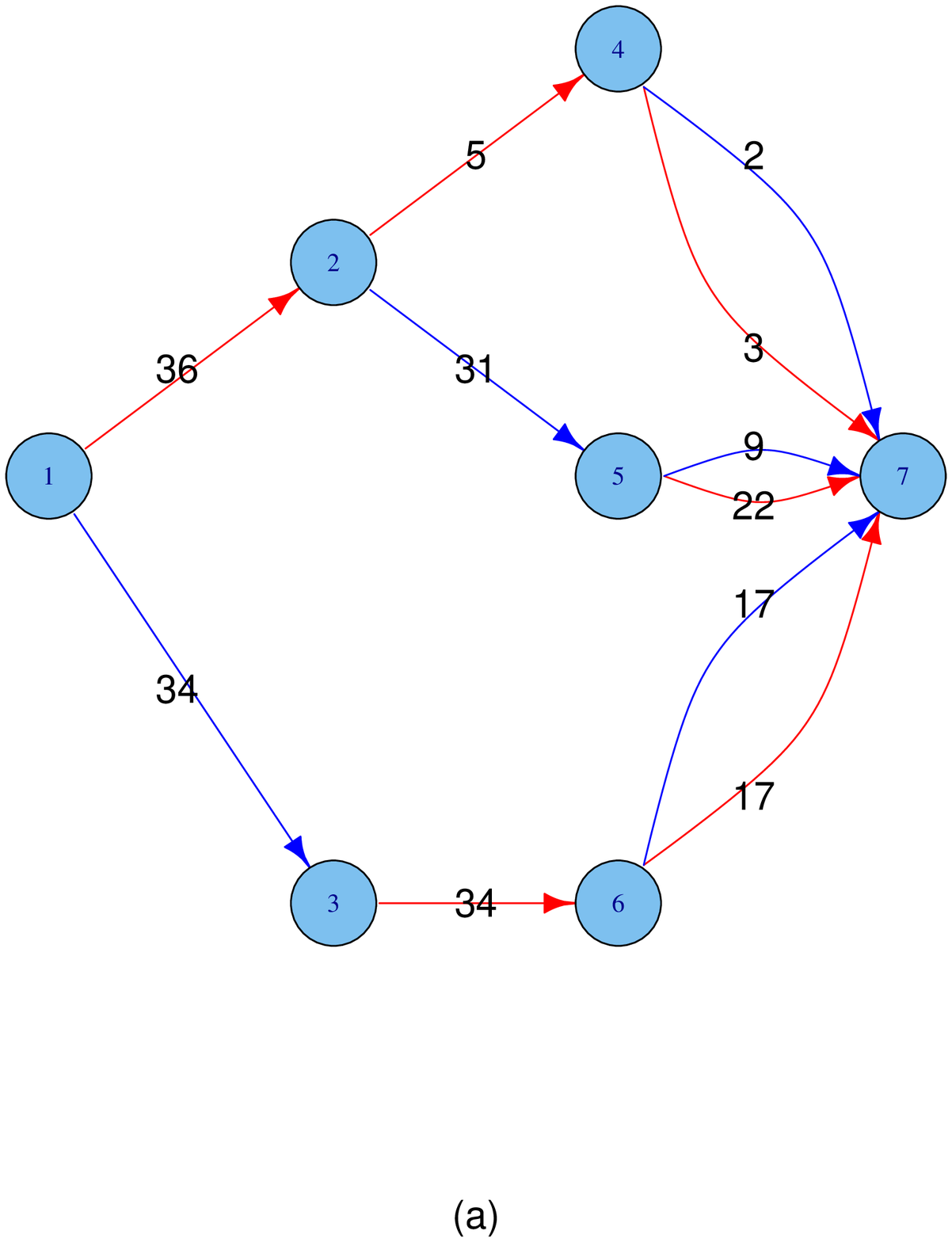}
\includegraphics[trim=2cm 0cm 0cm 2cm, clip=true, width=2.5in]{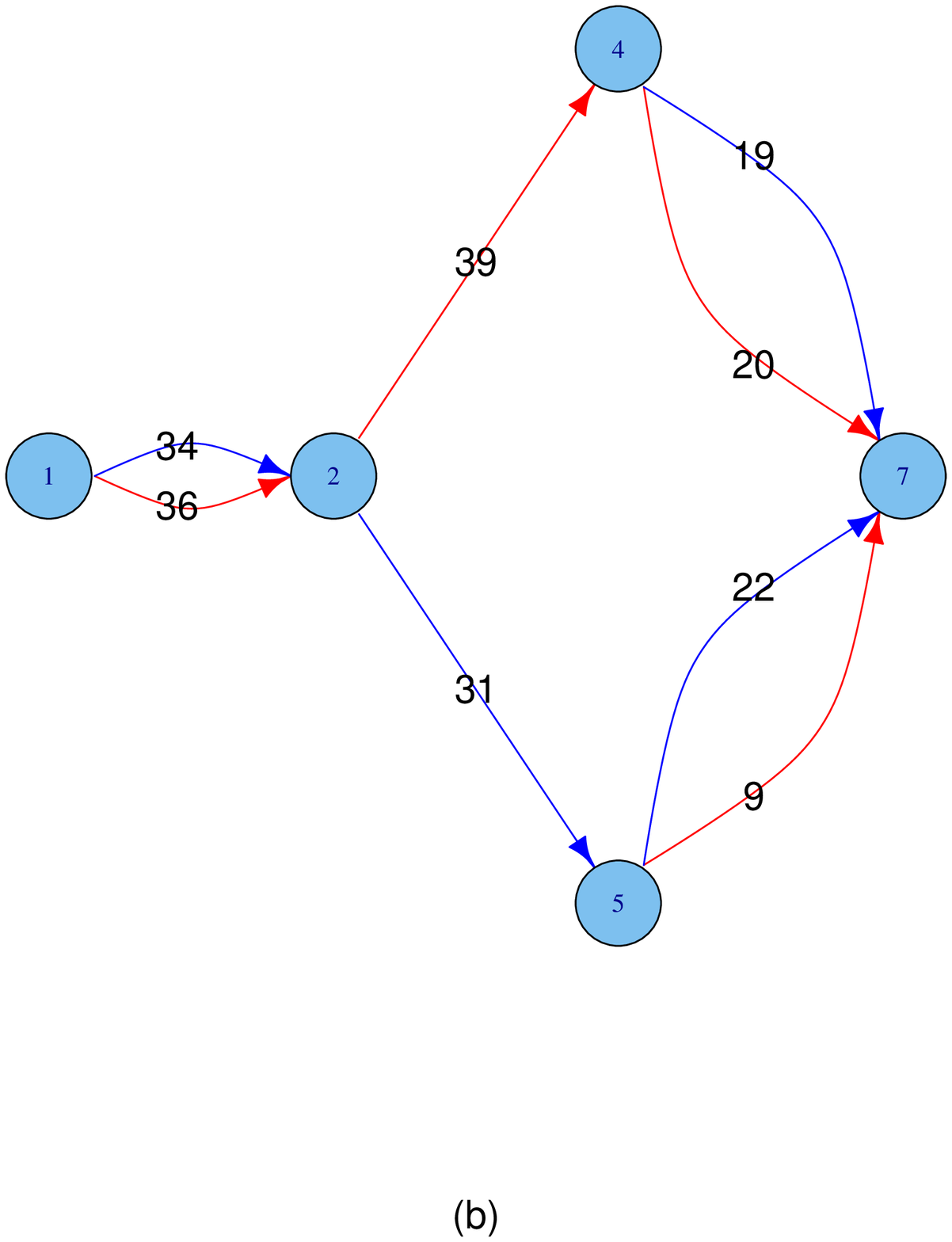}
\caption{State merging in an incomplete APFA. (b) is obtained from (a) by merging node 3 with node 2, and hence also node 6 with node 4.}
\label{fig:underspec}
\end{figure}

Now let $\cA$ and $\cA_0$ be the APFA shown in Figure~\ref{fig:underspec}. $\cA_0$ is obtained from $\cA$ by merging states 2 and 3. Suppose we wish to test whether these states can be merged. Note that $\cA$ is incomplete, and $\cA_0$ is not a submodel of $\cA$, but as we saw in Section~\ref{sec:merging} the underlying model spaces $(\cA_0)^+$ and $\cA^+$ are nested.

Using the approach just described, we can test for independence in the $2 \times 4$ contingency table shown in Table~\ref{ctab2}. This gives $G^2=67.288$ on three degrees of freedom.

\begin{table}
\caption{\label{ctab2}Comparing the distribution of $(X_2,X_3)$ given nodes 2 and 3: data in Figure~\ref{fig:underspec}(a).}
\centering
\fbox{%
\begin{tabular}{c|cccc}
  \hline
  source & (1,1) & (1,2) & (2,1) & (2,2) \\ \hline
  2 & 3 & 2 & 9 & 22 \\
  3 & 17 & 17 & 0 & 0 \\
  \hline
\end{tabular}}
\end{table}

The test statistic can be decomposed into a sum of $G^2$ statistics for two $ 2 \times 2$ tables, corresponding to the two state merging operations: that of state 3 with 2, and state 6 with 4. The tables are shown in Table~\ref{ctab3}. The first of these is formed from the counts on the outgoing edges from nodes 2 and 3 in Figure~\ref{fig:underspec}(a). The independence hypothesis states that the probabilities on the outgoing edges from state 2 are equal to those outgoing from state 3. The second table is similarly formed from the counts on the outgoing edges from nodes 4 and 6. We call these \textit{node-symbol} tables, since the rows correspond to the nodes being merged, and the columns to the symbols on the outgoing edges, and we call the $G^2$ tests for the two tables \textit{local LRTs}, since they only involve transitions from one level to the next.

Applying (\ref{eq:G2}) and (\ref{eq:df}) to the two tables results in test statistics of $67.1125$ and $0.1757$, each with one degree of freedom. Thus the test gives $G^2=67.288$ on two degrees of freedom. This differs from the previous test in that there are now two instead of three degrees of freedom. The decomposition implicitly builds on a re-parametrisation in which the parameter $\Pr(X_3=1|X_1=2, X_2=2)$  is inestimable, hence reducing the degrees of freedom by one. We call the quantity calculated in this way the \emph{adjusted} degrees of freedom. In larger APFA the adjusted and unadjusted degrees of freedom can differ substantially. The former quantity is preferable since it takes account of inestimability.

More generally, a likelihood ratio test of $\cA_0$ versus $\cA$, where $\cA_0$ is obtained by merging a set $s$ of nodes at level $i$ in $\cA$, may be computed as the sum of the local LRTs corresponding to the elements of $\cL(s)$. To see this we first re-express (\ref{eq:loglike}) as
\begin{equation}
\hat{\ell}(\cA)=\sum_{v \in V(\cA)} \sum_{e: s(e)=v} n(e) \log \hat{\pi}(e).
\label{eq:loglike1}
\end{equation}
\noindent
Note that for any $v \in V(\cA)$ that is not merged in $\cA_0$ and so not contained in any element of $\cL(s)$, the counts of the outgoing edges remain unchanged after merging and hence also the contribution $\sum_{e: s(e)=v} n(e) \log \hat{\pi}(e)$ to the log-likelihood remains unchanged. So it is sufficient to consider the node sets that are merged in $\cA_0$: for such a set, the contribution to the deviance (\ref{eq:deviance}) due to merging is equal to the likelihood ratio test for corresponding node-symbol table, as given in (\ref{eq:G2}), with degrees of freedom given by (\ref{eq:df}). Note also that if $\cL_{\cA}(s)$ and  $\cL_{\cA^+}(s)$ are the merge-lists computed in $\cA$ and $\cA^+$, then $\cL_{\cA^+}(s)$ will contain the node-pairs in $\cL_{\cA}(s)$ as well as node-pairs in which one or both nodes have zero node counts. These latter will not contribute to the log-likelihood or degrees of freedom, confirming that the computations can be based on $\cA$ alone.
\begin{table}
\caption{\label{ctab3}Decomposition of $G^2$ into two $2 \times 2$ tables.}
\centering
\fbox{%
\begin{tabular}{ccrc}
  \hline
   element of & $2 \times 2$ & $G^2$ & df \\
   $\cL(2,3)$  & table \\  \hline
  (2,3) & $\begin{array}{cc} 5 & 31 \\ 34 & 0 \end{array}$ &  67.112 & 1\\ \hline
   (4,6) & $\begin{array}{cc} 3 & 2 \\ 17 & 17 \end{array}$ & 0.176 & 1  \\ \hline
  sum  & & 67.288 & 2  \\
\end{tabular}}
\end{table}

When $|s|=2$, an alternative way to compute the test is to calculate $G^2$ using (\ref{eq:deviance}), and the adjusted degrees of freedom from $\cA_0$ as the sum of (outdegree\footnote{The outdegree of a node is the number of edges outgoing from the node.} minus one) over the nodes resulting from the merges in $\cL(s)$.

Up to now in this section we have described likelihood ratio tests for testing $\cA_0$ versus $\cA$, where $\cA_0$ is the submodel of $\cA$ formed by merging two nodes of $\cA$. More generally let $\cA_0$ be any submodel of $\cA$. Then for each level $i=1, \dots, p-1$, the level $i$ nodes of $\cA_0$ correspond to elements of a partition of the level $i$ nodes of $\cA$, in that such node is the result of merging the nodes in the corresponding element of the partition. It follows that the likelihood ratio test of $\cA_0$ versus $\cA$ can be decomposed into the sum over all levels and partition elements of the corresponding local $G^2$ quantities.

\section{Model Selection}
\label{sec:modelselect}
Suppose now that we have obtained independent samples $\bx^{(k)}=(x^{(k)}_1, \dots x^{(k)}_p)$ for $k=1 \dots N$ drawn from some
unknown APFA $\cA$, and we want to estimate (or select) $\cA$. \cite{Ron1998} describe a simple and efficient algorithm to do this.
It starts from the sample APFA, which is then simplified in a series of state merging operations. The intention is to merge two nodes $v$ and $w$ at level $i$ whenever (\ref{eq:condinda}) holds after merging, which implies that
\begin{equation}
\bX_{>i} \cip I_i | I_i \in \{v,w\}.
\label{eq:merging}
\end{equation}
\noindent
To assess this, a similarity score $\delta(v,w)$ between nodes $v$ and $w$, and a fixed threshold, $\mu$, are used. A small value of $\delta(v,w)$ means that the conditional distributions of $\bX_{>i}$ given $I_i=v$ and $I_i=w$ are similar. In particular, $v$ and $w$ are called \emph{similar} if $\delta(v,w) < \mu$: otherwise they are called \emph{dissimilar}. Dissimilar nodes are not merged.

The algorithm proceeds from levels 1 to $p-1$. At each level, the most similar pair of nodes is merged, and this is repeated until all the resulting nodes at the level are pairwise dissimilar. The algorithm then proceeds to the next level.

In \cite{Ron1998} the similarity score is defined as the maximum absolute value of the conditional probability differences for corresponding descendent nodes of $v$ and $w$, that is, nodes $x$ and $y$ for which there exist paths $v \rightarrow x$ and $w \rightarrow y$ with the same symbol sequence. Beagle uses the same similarity score but the threshold is allowed to vary, depending on the node counts $n(v)$ and $n(w)$ \citep{browning2007a}.

Here we sketch a natural alternative approach that is studied in more depth in \cite{Ankinakatte2013}. This is based on a penalized likelihood criterion
\begin{equation}
IC(\cA) = -2\hat{\ell}(\cA) + \alpha \dim(\cA)
\label{eq:IC}
\end{equation}
where  $\dim(\cA)$ is the number of free parameters under $\cA$, and $\alpha$ is a tuning parameter. For example, choosing $\alpha=2$ gives the Akaike information criterion \citep{akaike1974}, and choosing $\alpha=\log(N)$ gives the Bayesian information criterion \citep{schwarz1978}.  We define the \emph{penalized likelihood similarity score} as
\begin{eqnarray}
\nonumber
\delta_{pl}(v,w) & = & IC(\cA_0) - IC(\cA) \\
                       & = &  G^2 - \alpha k
\label{eq:g2sim}
\end{eqnarray}
\noindent
where $\cA_0$ is the APFA obtained after merging $v$ and $w$ in $\cA$, and $G^2$ and $k$ are the corresponding deviance statistic and adjusted degrees of freedom. The threshold is set to zero, so that two nodes are judged to be dissimilar whenever merging them would increase (\ref{eq:IC}).  Thus the selection algorithm using (\ref{eq:g2sim}) seeks to minimize (\ref{eq:IC}).

In \cite{Ankinakatte2013} the performance of this algorithm is compared to the algorithm implemented in Beagle in terms of both rate of convergence to the true model as $N \rightarrow\infty$ and prediction accuracy. The algorithm based on (\ref{eq:g2sim}) performs as well or better than that in Beagle in both respects.
\section{Conditional APFA models}
\label{sec: conditionalAPFA}
In this section we sketch how the framework may be extended to incorporate covariate information. For ease of exposition we assume that one covariate $\bz=(z^{(1)} \dots z^{(N)})$ is available, in addition to the $p$ discrete observed variables $\bx^{(k)}=(x^{(k)}_1, \dots x^{(k)}_p)$ for $k=1 \dots N$. We assume an APFA $\cA$ but allow the edge probabilities to depend on $\bz$, that is, by replacing $\pi(e)$ by $\pi(e \cd \bz)$, and adopting suitable parametric models for these conditional probabilities. Thus the likelihood of the sample, instead of
(\ref{eq: likelihood}), becomes
\begin{equation}
 \prod_{k=1 \dots N} \Pr(\bx^{(k)}\cd \theta, z^{(k)}) = \prod_{k=1 \dots N} \prod_{i = 1 \dots p} \pi(e_i^{(k)} \cd \theta, z^{(k)})
 \label{eq: condlikelihood}
\end{equation}
where $\theta$ is the parameter vector, and as before  $ \be^{(k)}= \sigma^{-1}(\bx^{(k)})$. The context-dependent conditional independence relations implied by $\cA$ involve conditioning on $Z$, that is, (\ref{eq:condinda}) becomes
\begin{equation}
\bX_{>i} \cip \bX_{\leq i} | I_i=w, Z=z
\label{eq:condinda2}
\end{equation}
for each $z$.

As a simple example, suppose that the covariate $Z$ is binary. Then the conditional model states that the data for the two groups are generated by the same APFA but with distinct sets of edge probabilities, $\pi(e \cd z=1)$ and $\pi(e \cd z=2)$.  The maximum likelihood estimates of these are the within-group relative frequencies of the corresponding edge counts. Using standard contingency table methods we can construct local LRTs for hypotheses of the type
\[
  X_{i+1} \cip I_i \cd I_i \in (v,w), Z
\]
for two level $i$ nodes, $v$ and $w$. This replaces the test of independence of $X_{i+1}$ and $I_i$ (given $ I_i \in (v,w)$) described in Section~\ref{sec:testing} with a test of conditional independence $X_{i+1}$ and $I_i$ given $Z$ (again, also given $ I_i \in (v,w)$). The test statistic and associated degrees of freedom are simply the sum of the corresponding within-group quantities (\ref{eq:G2}-\ref{eq:df}). To test $\cA_0$ versus $\cA$, where $\cA_0$ is formed by merging two level $i$ nodes $v$ and $w$, we consider the hypothesis
\[
  X_{> i} \cip I_i \cd I_i \in (v,w), Z.
\]
By the same logic as in Section~\ref{sec:testing}, the likelihood ratio test for this can be computed as the sum over $\cL(\{v,w\})$ of the corresponding local LRTs. Similarly, the incremental change in information criteria can be computed and used as the basis for the model selection algorithm of Section~\ref{sec:modelselect}.  It is well-known that marginal independence neither entails or is entailed by conditional independence \citep[Section 1.4]{edwards2000}, so the selected model may be simpler or more complex than that found using the unconditional approach.

The conditional tests and model selection process just described can be formulated in an alternative APFA framework in which the covariate $Z$ is included as a variable preceding  $X_1, \dots X_p$ (see Figure~\ref{fig:sampletree1} below). We omit the details. A comparison with the current approach would be valuable but is not attempted here.

Suppose now that $Z$ is continuous and the variables $X_1, \dots, X_p$ are binary. One choice of model for $\pi(e \cd z)$ is the logistic regression model
\begin{equation}
 \ln  {{ \pi(e \cd z) } \over {1 - \pi(e \cd z) }}  = a_e + b_e z
\label{eq: logist}
\end{equation}
where $a_e$ and $b_e$ are scalar parameters. (This would apply to one outedge $e$ of each node, say corresponding to $x_i=1$; for the other out-edge, say $\tilde{e}$, $\pi(\tilde{e} \cd z) = 1-\pi(e \cd z)$). To obtain maximum likelihood estimates of $a_e$ and $b_e$, a standard logistic regression algorithm is applied to the observations with $I_i=v$, where $v=s(e)$ is the source node of $e$.

Consider two level $i$ nodes, $v$ and $w$, with $e$ and $f$ being corresponding out-edges (for example, corresponding to $x_{i+1}=1$). We can construct a local LRT for the hypothesis
$
  X_{i+1} \cip I_i \cd I_i \in (v,w), Z
$
by applying logistic regression models to the observations with $I_i \in \{v,w\}$. Under the alternative the intercept and regression coefficients in (\ref{eq: logist}) for $I_i=v$ differ freely from those for $I_i=w$: under the null they are equal.  Again, a likelihood ratio test of $\cA_0$ versus $\cA$, where $\cA_0$ is obtained by merging $v$ and $w$, can be computed as the sum over $\cL(\{v,w\})$ of the corresponding local LRTs, and from this the incremental change in information criteria can be derived and used in the model selection algorithm of Section~\ref{sec:modelselect}.

Additional modelling possibilities and complexities may arise in the conditional setting. For example, it will often be useful to characterize the effects of the covariates in more detail.
We can examine for the individual edges $e$ in the binary case whether $\pi(e \cd z=1)= \pi(e \cd z=2)$, or in the continuous case whether $b_e=0$; with multiple covariates some kind of covariate selection procedure at each node could be used. Merging nodes requires that the associated conditional models have the same structure, so it is natural to select covariates after the graph is determined. We remark in passing that time-dependent covariates may be used in the same way, provided they are exogenous to the system.

\section{An Application}
\label{sec:application}
As part of the Six Cities study, a longitudinal study of the respiratory health effects of air pollutants \citep{ware1984}, the presence of absence of wheezing were recorded annually for a sample of 537 children from Steubenville, Ohio. The four variables, here denoted $W_1, \dots, W_{4}$, record the presence or absence of wheezing at ages 7, 8, 9 and 10. In addition, maternal smoking was recorded, here denoted $Z$, categorized as 1 if the mother smoked regularly and 0 otherwise. Although maternal smoking is a time-varying covariate, it is treated as fixed at its value in the first year of the study. In this section we apply the methods described above to these data. We first model the wheezing variables to gain insight into their dependence structure, then examine the possible effect of maternal smoking on the wheezing of their children.  See \cite{Ekholm1995} and \cite{fitzmaurice1993} for alternative analyses of these data.

The data are shown in Figure~\ref{fig:sampletree}. If we apply the algorithm of Section~\ref{sec:modelselect}, setting
$\alpha=\log(537)$ so as to minimize the Bayesian information criterion, we obtain the APFA shown in Figure~\ref{fig:wheeze}. The constraints to the joint distribution are shown in (\ref{eq:condp1}) and (\ref{eq:condp2}). At ages 9 and 10, the children fall into three groups: those with wheezing absent at each previous age, those with wheezing present at each previous age, and an intermediate group for which wheezing has sometimes been absent and sometimes present. Curiously, at both ages 9 and 10, the probability of wheezing is greater for the intermediate group than for the second group. It is also notable that the transition probabilities for those in which wheezing is absent, or present, at all ages increase monotonically over time, suggesting perhaps that subgroups of never- and always-wheezers are crystallizing out.

\begin{figure}[!th]
\begin{center}
\includegraphics[trim=2cm 6cm 0cm 4cm, clip=true, width=4.5in]{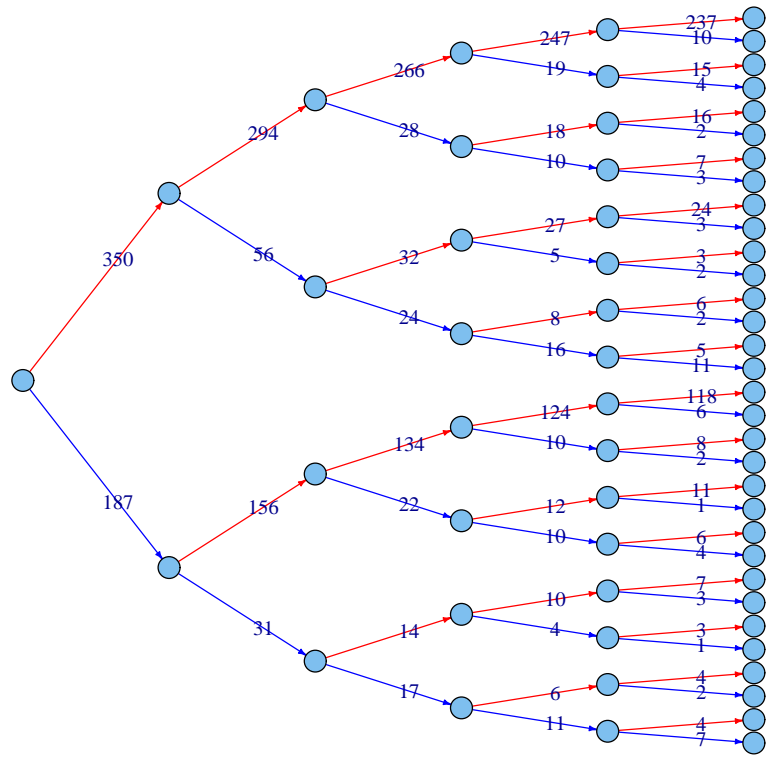}
\caption{A modified sample tree for the wheezing data, including the maternal smoking status at the first level, coded as non-smoker (red)
or smoker (blue). For the remaining edges, blue and red represent the presence and absence of wheezing. }
\label{fig:sampletree1}
\end{center}
\end{figure}

\begin{figure}[!ht]
\centering
\includegraphics[trim=0cm 0cm 0cm 2cm, clip=true, width=2.5in]{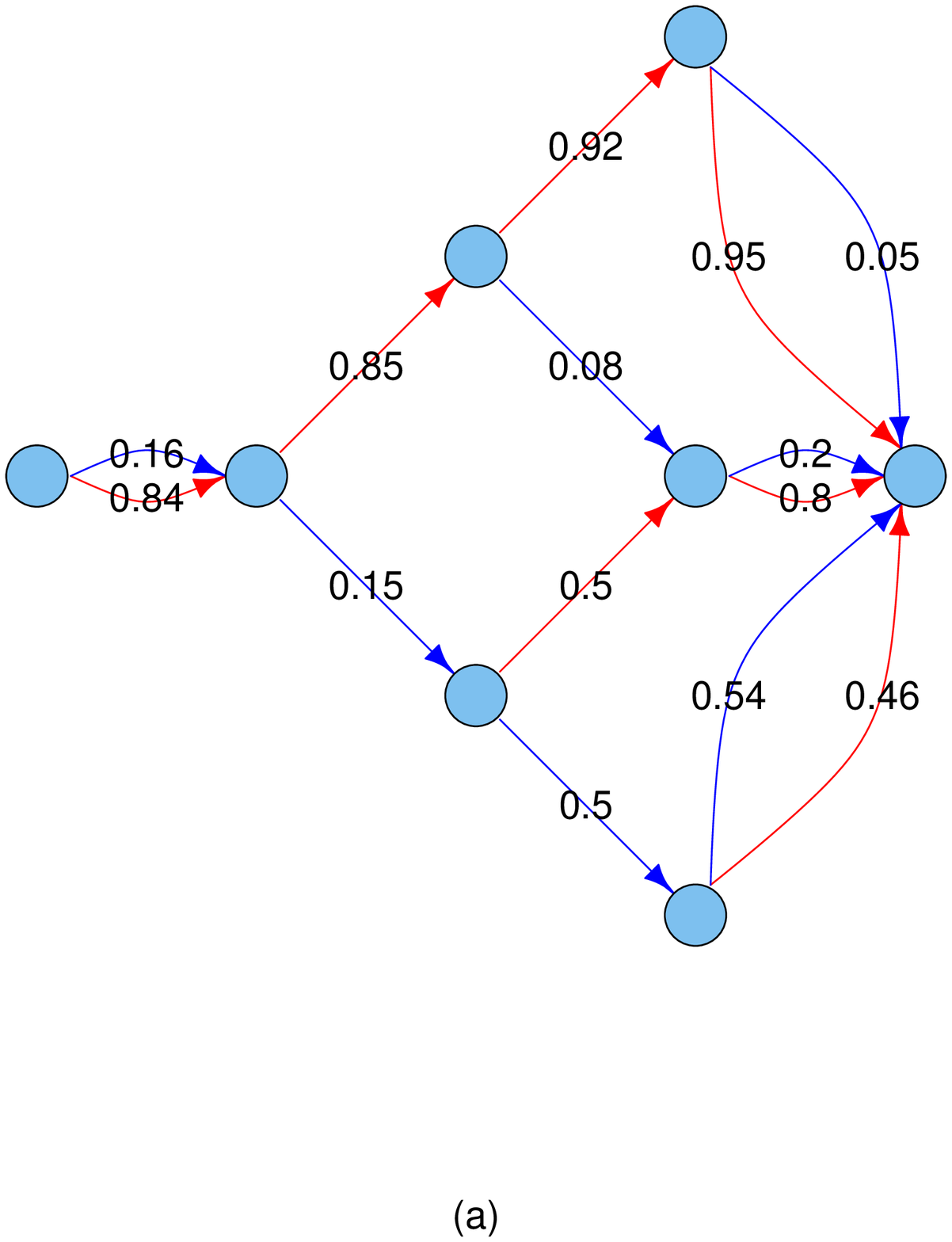}
\includegraphics[trim=0cm 0cm 0cm 2cm, clip=true, width=2.5in]{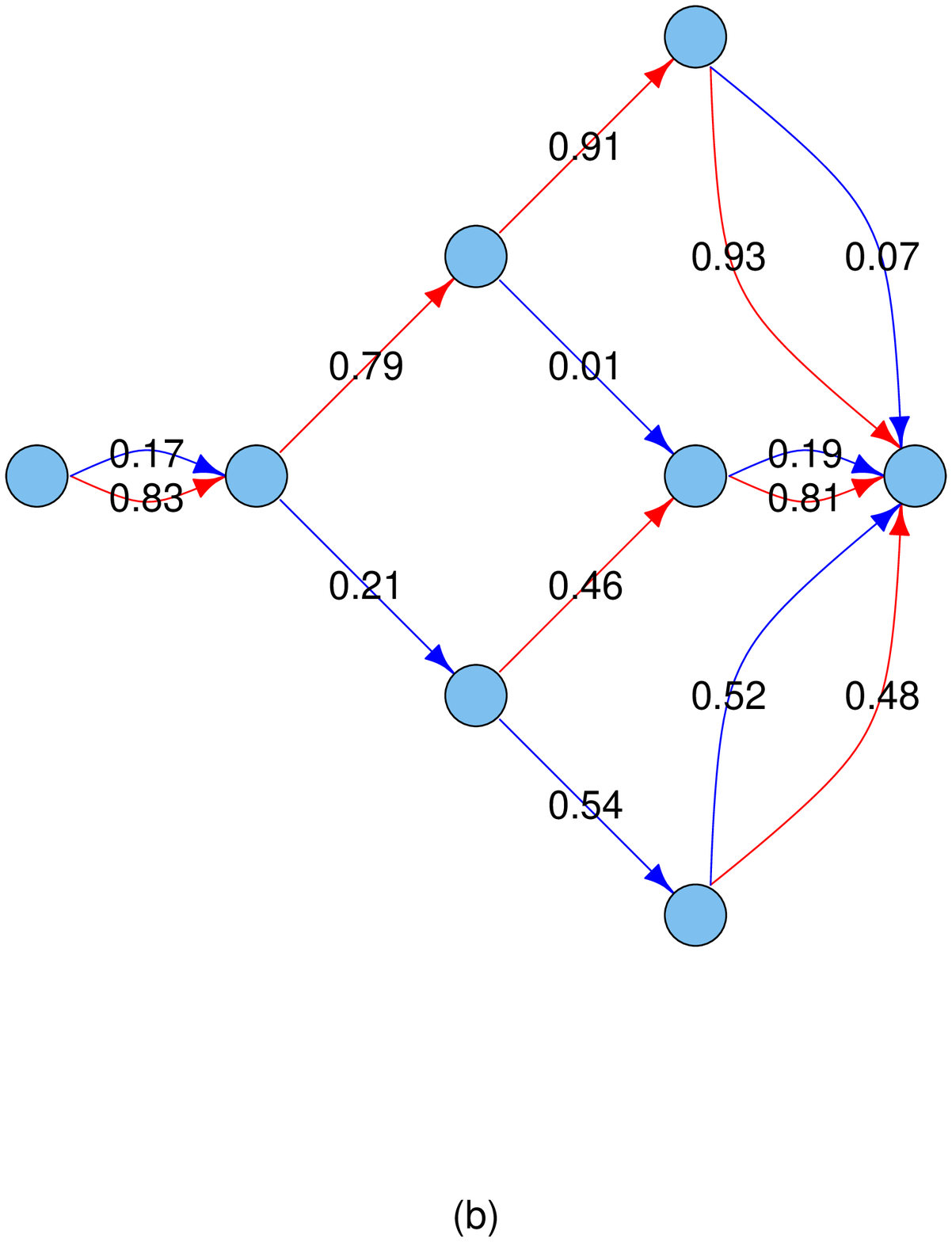}
\caption{A conditional APFA for the wheezing data: (a) and (b) show the transition probabilities for the children with mothers that were smokers, and non-smokers, respectively. As before, the blue and red edges represent the presence and absence of wheezing. }
\label{fig:wheeze_cova}
\end{figure}

To examine whether or not the maternal smoking affects the wheezing of her child, we first look at the modified sample tree shown in Figure~\ref{fig:sampletree1}. It is seen that there are 187 mothers that smoked regularly during the first year of the study and 350 that did not. We therefore assume a conditional APFA model, as described in Section~\ref{sec: conditionalAPFA}, and again apply the minimum BIC algorithm of Section~\ref{sec:modelselect}. This results in the model shown in Figure~\ref{fig:wheeze_cova}. This implies that
$W_1 \cip (W_2, W_3, W_4) \cd Z$, whereas under Figure~\ref{fig:simpleAPFA},  $W_1 \notcip (W_2, W_3, W_4)$: that is, conditional independence but marginal dependence. The transition probabilities for the children of maternal smokers and maternal non-smokers are very similar, suggesting that the conditional independence found is due to a power reduction rather than a significant covariate effect.

As a global test of equality of the conditional distributions of $W_1, \dots, W_{4}$ given $Z$ we can calculate the test for merging the two nodes at level one in the APFA formed from Figure~\ref{fig:sampletree1}: this gives $G^2 = 7.36$ with 15 degrees of freedom. So there is little evidence of any effect of maternal smoking on the wheezing of the children.

We remark that the conditional model could have used the annual recordings of maternal smoking, that is, as a time-varying covariate. However, this would require assuming that the maternal smoking was not affected by their children's wheezing, which would seem uncharitable.

\section{Related Markov models}
\label{sec:modelclasses}
As we saw above, APFA are context-specific graphical models for discrete longitudinal data. There has been substantial recent interest in extending Markov and Bayesian network modelling techniques so as to incorporate context-specific information. These allow independence structure to vary locally in ways that are not captured by standard Markov or Bayesian networks. Examples include \cite{boutilier1996, corander2003, eriksen1999, eriksen2005,hara2012, hojsgaard2003, hojsgaard2004, jaeger2004, myers2007}.

In this section we first briefly relate APFA to the transition models often used in the analysis of discrete longitudinal data. Then we compare APFA to discrete graphical models, both undirected (Markov networks) and directed (Bayesian networks). Finally we briefly relate APFA to chain event graphs. One motivation for these comparisons is to understand better the strengths and limitations of the different types of model, so as to inform choice between these in specific applications.

\begin{figure}[!ht]
\centering
\includegraphics[trim=1cm 1cm 0cm 1cm, clip=true, width=1.5in]{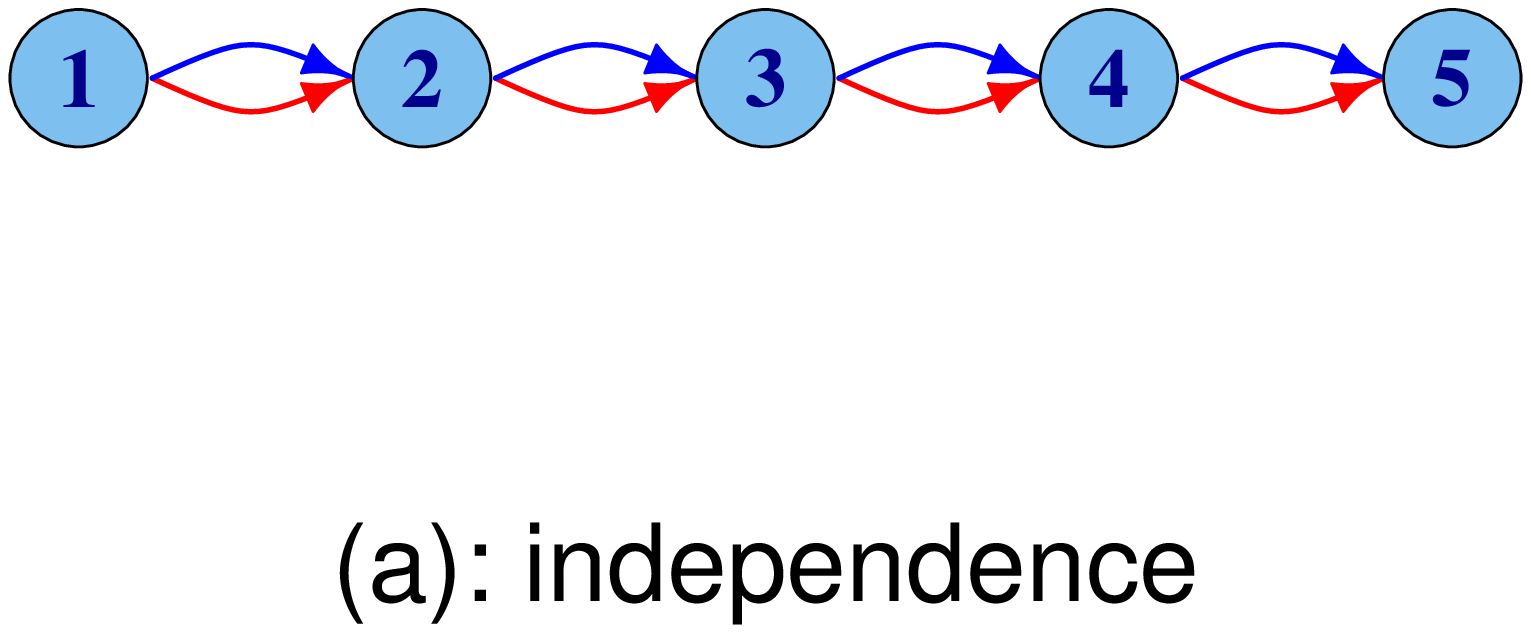}
\includegraphics[trim=1cm 1cm 0cm 1cm, clip=true, width=1.5in]{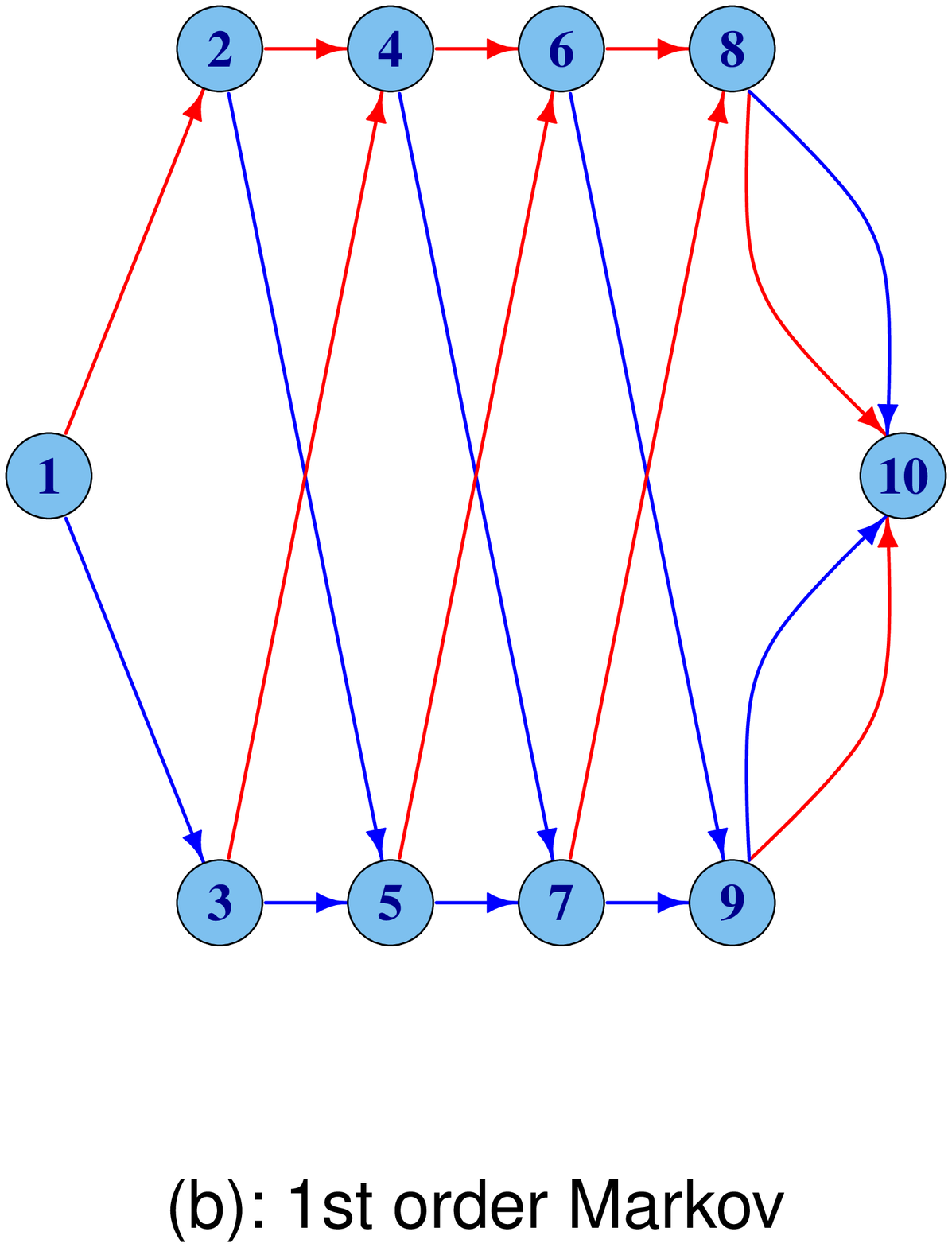}
\includegraphics[trim=1cm 1cm 0cm 1cm, clip=true, width=1.5in]{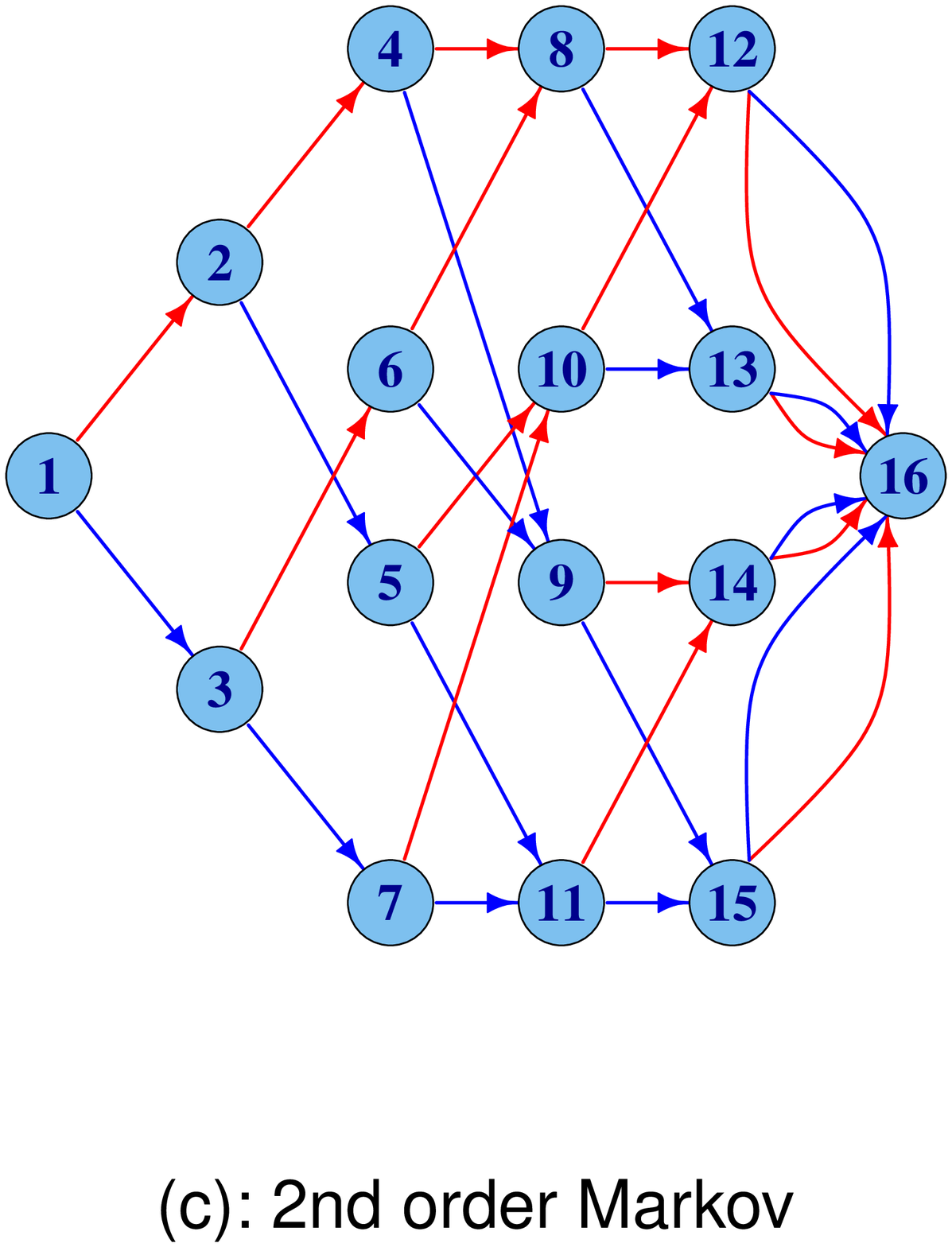}
\includegraphics[trim=1cm 1cm 0cm 1cm, clip=true, width=1.5in]{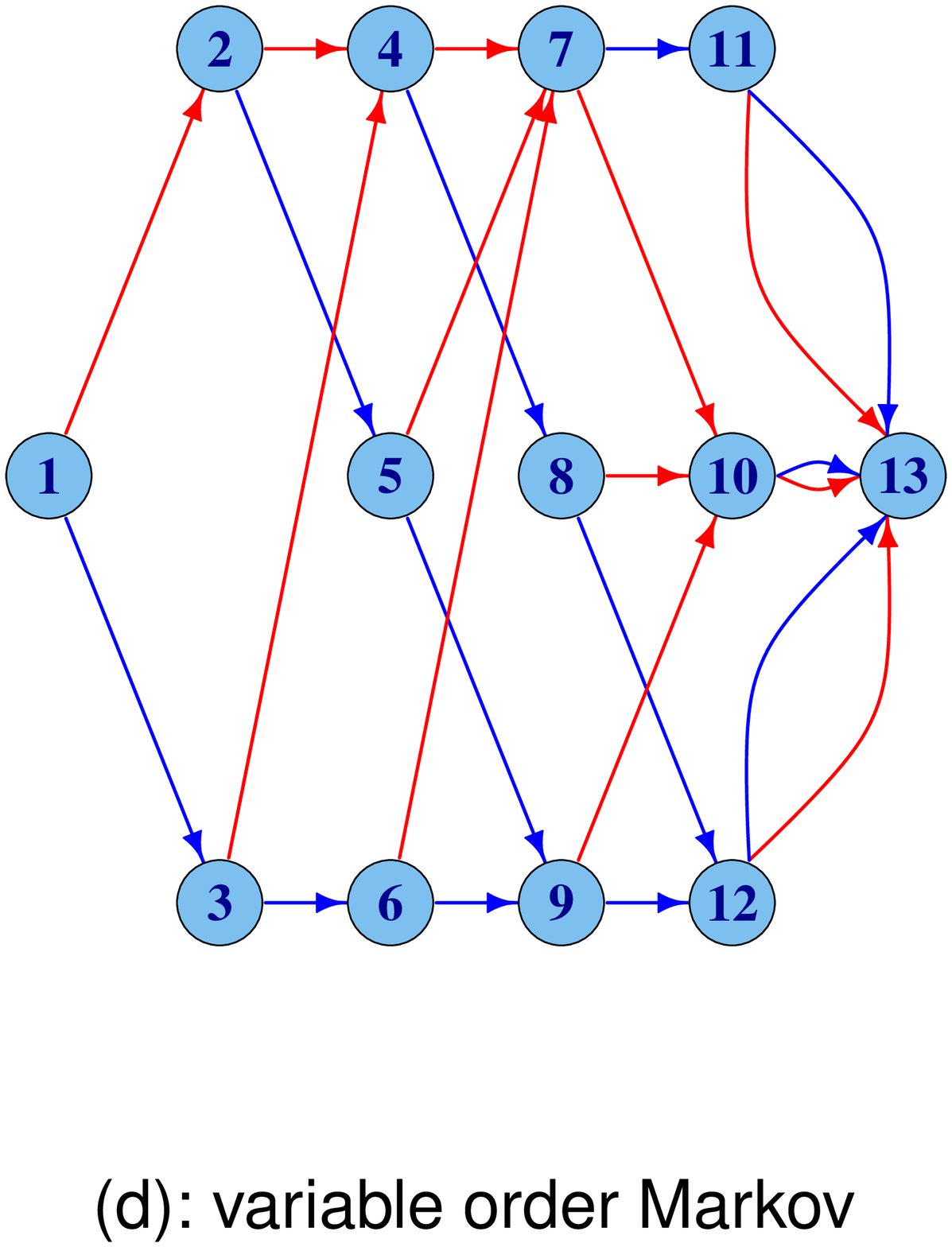}
\includegraphics[trim=1cm 1cm 0cm 1cm, clip=true, width=1.5in]{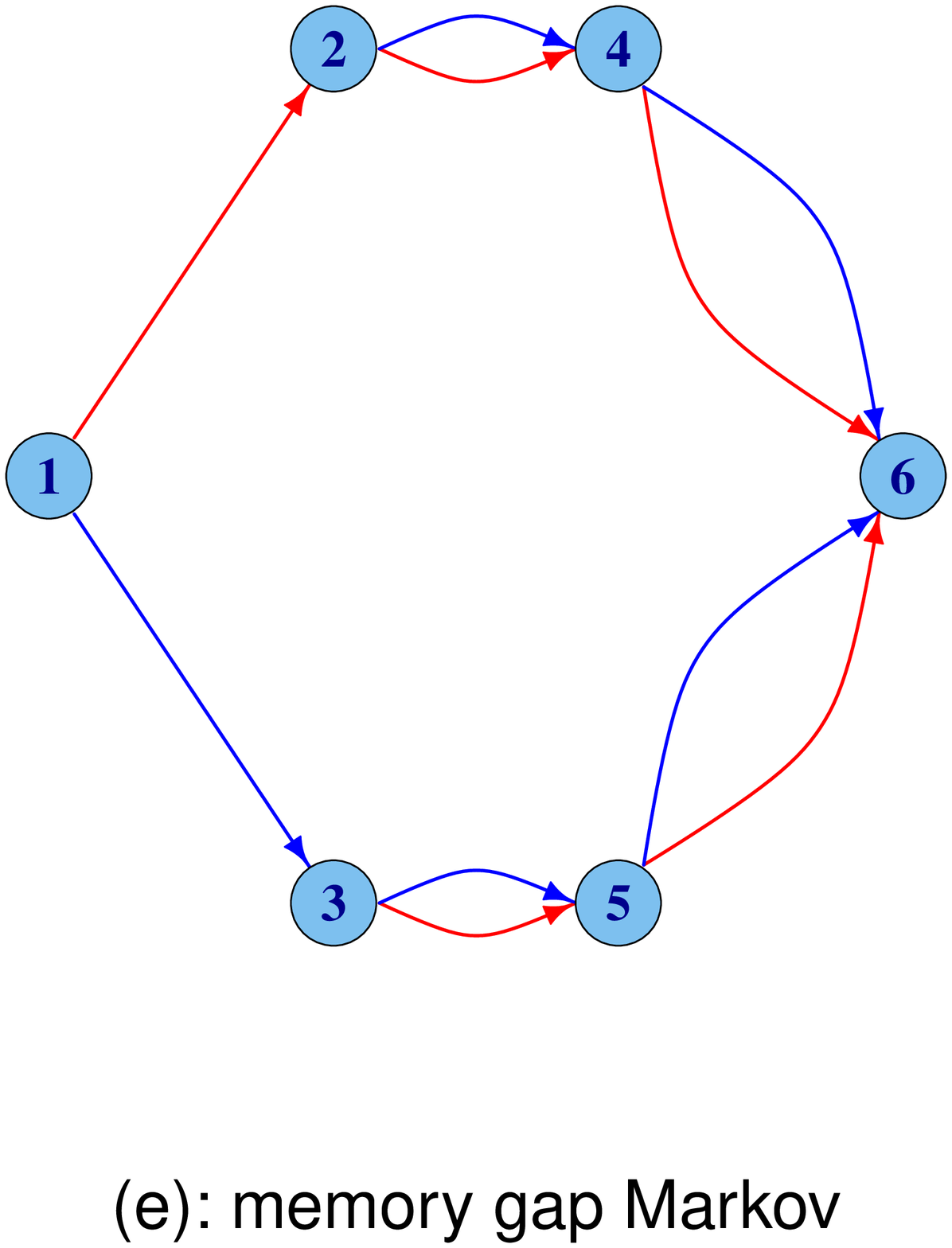}
\caption{Some APFA that are equivalent to standard Markov models.}
\label{fig:gMex}
\end{figure}

Transition models focus on the conditional distribution of the response $X_{j}$ at time-point $j$ given the prior responses $X_1, \dots, X_{j-1}$ and possibly covariates, say $Z_1, \dots, Z_s$. Generally the conditional distribution involves only on the previous $q$ responses, where $q$ is called the \emph{order} of the model. The $q$ prior responses and $s$ covariates are treated on an equal footing as explanatory variables in a convenient parametric model for $X_{j}$.  There is much freedom in the choice of parametric model: for example, generalized linear models may be used \citep[Chapter 10]{diggle2013}. Stationarity is generally assumed, so that the conditional distributions are constant over the time interval spanned by the data. A simple example for binary data and no covariates is a $q$th order Markov chain, in which the transition probabilities are specified by a $2^q$ table of conditional probabilities. In contrast, APFA are non-stationary, requiring no assumption of constancy of conditional distributions over time. Furthermore APFA have no fixed order, but rather allow the length of dependence to vary, as in variable order Markov chains (see Figure~\ref{fig:gMex}(d) below). Thus APFA are appropriate for non-stationary data that exhibit long-range dependences.

We now turn to discrete graphical models. We are interested in equivalences between three model classes for $p$ discrete variables, $\mX= X_1, \dots X_p$: APFA models, which we denote $\Theta$; directed graphical models (Bayesian networks), which we denote $\Delta$; and undirected graphical models, which we denote $\Upsilon$. The models in $\Delta$ and $\Upsilon$ have $p$ nodes, corresponding to $X_1, \dots X_p$, and which we label $1, \dots p$. For models in $\Delta$ we require that the directions are consistent with the variable ordering, that is, there may exist a directed edge from node $i$ to $j$ only when $i < j$.

We now examine equivalences between $\Theta$, $\Delta$ and $\Upsilon$ in more detail.

Directed and undirected graphical models are characterized by sets of conditional independence constraints of the form
\begin{equation}
\bX_B \cip \bX_C \cd \bX_A
\label{eq:UGcondind}
\end{equation}
for certain set triplets $(A,B,C)$. Such constraints may be re-written in less abbreviated form as
\begin{equation}
\bX_B \cip \bX_C \cd \bX_A=\bx_A
\label{eq:UGcondind1}
\end{equation}
for each $\bx_A \in \mX_A$, the sample space of $\bX_A$. In contrast, APFA are characterized by context-specific conditional independence constraints of the form (\ref{eq:condind}), that is
\begin{equation}
\bX_{> i} \cip \bX_{\leq i} \cd \bX_{\leq i} \in \cC(w)
\label{eq:condind1}
\end{equation}
for each node $w$ at level $i$. We recall that ${\cC(w)}=\{\sigma(\be): \be \in \cP(w)\}$, where $\cP(w)$ is the set of paths from the root to $w$.
Two special cases of (\ref{eq:condind1}) should be noted. Firstly, when $\cP(w)$ contains only one path, (\ref{eq:condind1}) is devoid of content. Secondly, when $w$ is the only node at level $i$, the event $\bX_{\leq i} \in \cC(w)$ has probability one, so (\ref{eq:condind1}) states that
$\bX_{\leq i}$ and $\bX_{>i}$ are marginally independent. For example, in Figure~\ref{fig:gMex}(a), there is a single node at levels $1$ to $3$, so three marginal independence statements hold: $X_1 \cip (X_2,X_3,X_4)$, $(X_1,X_2) \cip (X_3,X_4)$ and $(X_1,X_2,X_3) \cip X_4$.

To relate (\ref{eq:condind1}) to (\ref{eq:UGcondind1}), let ${\cE}_i(\cA)$ be the set of paths in $\cA$ that start at the root and end at a node at level $i$. For a path $\be \in {\cE}_i(\cA)$ and $A \subset \{1, \dots i\}$, let $\sigma(\be)_A$ be the subvector of $\sigma(\be)$ corresponding to $A$. For an $\bx_A \in \mX_A(\cA)=\{\sigma(\be)_A: \be \in {\cE}_i(\cA)\}$, let $Q(\bx_A) = \{\be \in {\cE}_i(\cA): \sigma(\be)_A=\bx_A\}$. Whenever $\cP(w)= Q(\bx_A)$ for some $\bx_A\in \mX_A(\cA)$, the event $\bX_{\leq i} \in \cC(w)$ is equivalent to the event $\bX_A=\bx_A$ and so (\ref{eq:condind1}) is equivalent to a conditional independence statement of the form
(\ref{eq:UGcondind1}) with $B=\{1, \dots i\} \setminus A$, and $C=\{i+1, \dots ,p\}$. When this is true for each node at level $i$, a conditional independence statement of the form (\ref{eq:UGcondind}) with the same $B$ and $C$ holds. Thus for an APFA to be equivalent to a graphical model it is necessary that the following property holds: for each $i$ ($1 \leq i \leq p-1$) there is a set $A(i) \subseteq \{1, \dots i\}$, such that for each level $i$ node $w$, $\cP(w)= Q(\bx_{A(i)})$ for some $\bx_{A(i)} \in \mX_{A(i)}(\cA)$. We call this \emph{property Q}. See Figure~\ref{fig:gMex} (b), (c) and (e) for some examples.

Let $\cA$  be an APFA with property Q, where the sets $A(i)$ are chosen to be maximal sets for which the property holds. We now show that for each $i \in \{1, \dots, p-1\}$
\begin{equation}
A(i) \subseteq A(i-1) \cup \{i\}
\label{eqn:running}
\end{equation}
where $A(0)$ is taken to be the null set. Suppose that $A(i) \not \subseteq  A(i-1) \cup \{i\}$ and let $B = A(i) \setminus (A(i-1) \cup \{i\})$. Then $\bx_B$ is constant for all paths in $\cP(w)$ for each level $i$ node $w$, and hence also constant for all paths in $\cP(v)$ for each level $i-1$ node $v$, contradicting the maximality of $A(i-1)$. Thus (\ref{eqn:running}) holds as stated.

To construct an equivalent model $\cG \in \Delta$, set $\pa(i)$\footnote{The parents $\pa(v)$ of a node $v$ in a directed graph are the nodes $w$ for which there exists an edge from $w$ to $v$. The adjacency set $\adj(v)$ of a node $v$ in an undirected graph is the set of nodes $w$ for which there exists an edge between $v$ and $w$.} using $\pa(i)=A(i-1)$, for $i=2, \dots p$. Note the sets $A(i)$ are complete in $\cG$, since $j \not \rightarrow k$ with $j<k<i$ implies $j \not \in A(k)$ and hence from  (\ref{eqn:running}) that $j \not \in A(i)$.  So $\cG$ has no immoralities and is therefore Markov equivalent to the undirected graphical model with the same skeleton \citep{gMwR2012}.  This undirected graphical model is decomposable.

Conversely, given a model in $\Delta$ satisfying $\pa(i) \subseteq \pa(i-1) \cup \{i-1\}$ for $i=2, \dots p$, we can construct a model in $\Theta$ with property Q for $A(i)=\pa(i+1)$, for  $i=1, \dots p-1$. Clearly we can do this for level 1: assume that we have done it up to level $j$. Then we define $I_{j+1}$ as the partition of the set of all combinations of the values of $I_j$ and $X_{j+1}$ that corresponds to $A(j+1)$. Hence the result follows by induction.

We have shown the following results.

\textbf{Theorem 1.}\emph{
An APFA  in $\Theta$ is equivalent to a directed graphical model in $\Delta$ (or an undirected graphical model in $\Upsilon$) if and only if it has property Q.}

\textbf{Theorem 2.}\emph{
A directed graphical model in $\Delta$ is equivalent to an APFA in $\Theta$ if and only if it satisfies $\pa(i) \subseteq \pa(i-1) \cup \{i\}$ for $i=2, \dots p$.}

\textbf{Theorem 3.}\emph{
An undirected graphical model in $\Upsilon$ is equivalent to an APFA in $\Theta$ if and only if it satisfies
$\adj(i) \cap\{1, \dots i-1\} \subseteq (\adj(i-1)  \cap\{1, \dots i-2\}) \cup \{i\}$ for $i=2, \dots p$.}

Figure~\ref{fig:gMex} shows some examples.

As mentioned above,  in Figure~\ref{fig:gMex}(a), the above conditions apply with  $A(i)=\emptyset$ for $i=1,2,3$,
implying three marginal independence statements hold: $X_1 \cip (X_2,X_3,X_4)$, $(X_1,X_2) \cip (X_3,X_4)$ and $(X_1,X_2,X_3) \cip X_4$. These imply complete independence. This model is equivalent to the undirected graphical model with four isolated nodes.

For the APFA shown in Figure~\ref{fig:gMex}(b), the above conditions apply with $A(i)=\{i\}$ for $i=1,2,3$. These give rise to three
conditional independences: $X_1 \cip (X_3,X_4,X_5) \cd X_2$, $(X_1,X_2) \cip (X_4,X_5) \cd X_3$, and $(X_1,X_2,X_3) \cip X_5 \cd X_4$.
This model is equivalent to a first order Markov model.

Figure~\ref{fig:gMex}(c) corresponds similarly to a second-order Markov model, with $A(i)=\{i-1,i\}$ for $i=1 \dots 4$. For the nodes at the first two levels, $\cP(w)$ contains one path only, and hence (\ref{eq:condind1}) is empty. From levels 3 and 4 we obtain
$X_1 \cip (X_4,X_5) \cd (X_2,X_3)$ and $(X_1,X_2) \cip X_5 \cd (X_3,X_4)$.

Figure~\ref{fig:gMex}(d)  represents a variable length Markov chain (VLMC). For these each node $w$ at level $i$ has $\cP(w)= Q(\bx_A)$ but for possibly different $A$'s. Also, implicit in the concept is that for level $i$, the A sets take the form $i-k, \dots i$ for some $k>0$ representing the memory length. At levels $i=2,3,4$, the paths to one node are characterized by $X_i=1$ (red=1) to another by $(X_{i-1}, X_i)= (1,2)$ and to the third by $(X_{i-1}, X_i)= (2,2)$.
Under the model
\begin{equation}
  \Pr(X_{\geq i}|X_{<i}) =  \left\{
  \begin{array}{@{}ll@{}}
     \Pr(X_{\geq i}|X_{i-1}) & \mbox{if $X_{i-1}=1$} \\
     \Pr(X_{\geq i}|X_{i-1},X_{i-2}) &   \mbox{if $X_{i-1}=2$}
  \end{array}
  \right.
\end{equation}
for $i=2,3,4$. Note since it does not satisfy the stated conditions it is not equivalent to a undirected or directed Markov model. It is however a submodel of the second order Markov model shown in Figure~\ref{fig:gMex}(d), and is obtained by merging nodes 4 with 6, 8 with 10 and 12 with 14.

Figure~\ref{fig:gMex}(e) exemplifies what might be called a memory gap Markov chain. It satisfies the stated conditions with $A(1)=A(2)=\{1\}$, and gives rise to the conditional independence $X_3 \cip X_2 | X_1$.

As mentioned above, \cite{smith2008} recently and independently introduced a class of models for discrete longitudinal data called chain event graphs. These include APFA as a special case, but also allow edges between non-adjacent levels, and focus on applications in which structural zeroes occur. The models are intended to be elicited from domain experts rather than selected on the basis of data samples, but model selection approaches have also been described \citep{freeman2011, cowell2013, silander2013}. In \cite{thwaites2011} various forms of conditional independence relations that hold under the models are studied. Other aspects have been developed in a series of papers: \cite{thwaites2010, thwaites2012, barclay2012, barclay2012a, riccomagno2005, riccomagno2009}.

\section{Discussion}
The intention of this paper has been to describe APFA from a statistical perspective. This has led to several contributions which we believe to be novel, including the basic results on hypothesis testing, the modification of the algorithm of \citet{Ron1998} to minimize information criteria, conditional APFA models, and the characterization of equivalences with graphical models. In this section we discuss some more general aspects and issues.

APFA are appropriate for discrete longitudinal data that are non-stationary and exhibit long-range dependences. They assume that the variables are measured at common times (or, in the case of genomic data, at common spatial positions). They scale well to high-dimensional data, since the model selection algorithm can be implemented very efficiently. We can report, for example, that when applied to a animal genetics data set with $N=16310$ observations of $p=44991$ binary variables, Beagle took just over 10 minutes of computing time to select an APFA.

As mentioned in Sections~\ref{sec:apfa} and \ref{sec:merging}, we have assumed that when combinations of variable values are not present in the data this is due to random rather than structural zeroes. The issue arises because the model selection procedure often involves comparisons between incomplete APFA that are ostensibly non-nested. Assuming that zeroes are random rather than structural allows us to regard an incomplete APFA as a computationally convenient representation of a larger complete APFA, so that the comparisons are between nested models and fall within a standard statistical framework. There is an analogy with backward selection procedures in contingency table analysis: initial complex models often have zero fitted counts, but these are not regarded as structural zeroes and are ignored in the subsequent selection process.

Nevertheless, structural zeroes may well occur: how does this affect the proposed methods? As a simple example consider $2$ binary variables, $X_1$ and $X_2$, and suppose that when $X_1=2$, $X_2$ cannot be 2. In small samples, the selection procedure may choose a model in which $X_1$ and $X_2$ are independent, which is in conflict with the structural zero. A possible remedy would be to modify the model selection algorithm so that state merges that conflict with any of a set of prior known structural zeroes are disallowed. Alternatively, approaches building on the concept of quasi-independence \citep[Chapter 5]{bishop2007} might be considered.

A question not mentioned above is how to apply the methods in the presence of missing data. A crude approach which may sometimes be adequate is to treat missing values as extra levels of the categorical variables. In Beagle an iterative algorithm is implemented that solves a specific genetic problem, combining model selection, phase estimation and imputation \citep{browning2007a}: see also \cite{cawley2003}. Algorithms need to be developed that in a more general setting maximize the marginal likelihood for a given APFA, impute missing values for a given APFA, and select an APFA in the presence of missing data.

Other issues also deserve further study. For example, can the model selection algorithm described in Section~\ref{sec:modelselect} be modified to preferentially select graphical models, or to handle ordinal discrete variables? Is there a minimal graphical model containing a given APFA as a submodel and if so can it be identified efficiently? Are methods developed for simplifying PFA \citep{dupont2000, thollard2008} useful for APFA? Also issues of causal analysis and interpretation \citep{thwaites2010} deserve further study.

\section{Acknowledgements}
Conversations with Poul Svante Eriksen and suggestions from an anonymous referee are gratefully acknowledged.

\bibliography{APFA2}


\section*{Appendix A: proof of Assertion I}

Let  $\cA^+$ be the completion of an APFA $\cA$, as defined in Section~\ref{sec:merging}, and let $\cA^m$ be the result of merging states $v$ and $w$ at level $i<p$ in $\cA$. We wish to show that  $(\cA^m)^+$ is a submodel of  $\cA^+$. Since $\cA^+$ is complete, the merging operation only imposes additional constraints of the form (\ref{eq:condind}) and hence  $(\cA^+)^m$ is a submodel of  $\cA^+$. It remains to show that $(\cA^+)^m=  (\cA^m)^+$.

Write $\cA=(V,E)$, $\cA^m=(V^m, E^m)$ and $\cA^+=(V \cup V_0, E \cup E_0) $ where $V_0$ and $E_0$ are the sets of nodes and edges added to $\cA$ in the process of completion. Since $\cA^+$ is complete, for any path $P$ from $v$ to the sink in $\cA^+$ there exists a path $Q$ from $w$ to the sink with the same symbol sequence.
In $(\cA^+)^m$ all such path-pairs $(P,Q)$ are merged. Write $P$ and $Q$ as
\[
  (v_0) \ e_1 \ (v_1) \ e_2 \ (v_2) \ \dots  \ e_q  \ (v_q) \ \ \mbox{and} \ \
  (w_0) \ f_1 \ (w_1) \ f_2 \ (w_2) \ \dots \ f_q  \ (w_q)
\]
where $v=v_0, v_1, \dots v_q$ and $w=w_0, w_1, \dots w_q$ are states,  $e_1, \dots e_q$ and $f_1, \dots f_q$ are edges, $q=p-i$ and both $v_q$ and $w_q$ represent the sink. Let $j_v = \max\{k \in \{0, \dots q-1\} : v_k \in V\}$ and similarly $j_w = \max\{k \in \{0, \dots q-1\}: w_k \in V\}$.
Then $j \leq j_v \Leftrightarrow v_j \in V$  for $j<q$, and $j \leq j_v  \Leftrightarrow e_j \in E$ for $j \leq q$. Similarly,
 $j \leq j_w \Leftrightarrow w_j \in V$  for $j<q$, and $j \leq j_w  \Leftrightarrow f_j \in E$ for $j \leq q$.

Let $j^- = \min(j_v, j_w)$. The sub-paths $\{v_j: j \leq j^-\}$ and $\{w_j: j \leq j^-\}$ are in $\cA$ and contain precisely the
node- and edge pairs that are merged in $\cA^m$.  When $j_w < j_v$ the sub-paths $\{v_j: j_w <j< j_v\}$ has nodes and edges in $V$ and $E$ respectively,
and  $\{w_j: j_w < j < j_v)\}$ has nodes and edges in $V_0$ and $E_0$ respectively. The converse holds when $j_v < j_w$. Identify the resulting merged nodes and edges with the corresponding nodes and edges in $V$ and $E$. Then it follows that the subgraph of $(\cA^+)^m$ induced by $V$ is identical to $\cA^m$.

Let $j^+ = \max(j_v, j_w)$. If $j^+ < p$ then the sub-paths $\{v_j: j^+ < j < p\}$ and $\{w_j: j^+ < j < p \}$ are both in $V_0$ and $E_0$, and
$v_{j^+}$ and $v_{j^+}$ both are incomplete in $\cA$, without an outgoing edge with the symbol $\sigma(e_{j^+} +1)=\sigma(f_{j^+} +1)$. Hence the node resulting from the merge of $v_{j^+}$ and $w_{j^+}$ in $\cA^m$ is also incomplete. Call this node $z$ in $\cA^m$ and $(\cA^+)^m$. Since $\cA^+$ is complete, the descendent graphs of $z$ in $(\cA^+)^m$ are complete trees (with the last level contracted to the sink), and so are identical with the descendent graphs of $z$ in  $(\cA^m)^+$. Hence $(\cA^+)^m=  (\cA^m)^+$. 

\end{document}